\documentclass[11pt,amsfonts]{article}
\usepackage{graphicx}
\usepackage{latexsym}
\usepackage{amssymb}
\usepackage{amsmath}
\usepackage{layout}

\newtheorem{prop}{Proposition}[section]
\newtheorem{rem}[prop]{Remark}
\newtheorem{proof}{Proof}
\newtheorem{df}{Definition}[section]
\newtheorem{thm}[prop]{Theorem}
\newtheorem{lem}[prop]{Lemma}
\newtheorem{cor}[prop]{Corollary}
\newtheorem{ex}[prop]{Example}
\def\real{{\mathord{{\rm I\kern-2.8pt R}}}}        
\def\inte{{\mathord{{\rm I\kern-2.8pt N}}}}

\def\sZZ{{\rm Z\kern-2.8ptem{}Z}}

\def\z{{\mathchoice
  {\sZZ}
  {\sZZ}
  {\rm Z\kern-0.30em{}Z}
  {\rm Z\kern-0.25em{}Z} }}
\def\sQQ{{\kern 0.27em \vrule height1.45ex width0.03em depth0em
          \kern-0.30em \rm Q}}
\def\qu{{\mathchoice
    {\sQQ}
    {\sQQ}
  {\kern 0.225em \vrule height1.05ex width0.025em depth0em \kern-0.25em \rm Q}
  {\kern 0.180em \vrule height0.78ex width0.020em depth0em \kern-0.20em \rm Q}
        }}
\def\sCC{{\kern 0.27em \vrule height1.45ex width0.03em depth0em
          \kern-0.30em \rm C}}
\def\complex{{\mathchoice
    {\sCC}
    {\sCC}
  {\kern 0.225em \vrule height1.05ex width0.025em depth0em \kern-0.25em \rm C}
  {\kern 0.180em \vrule height0.78ex width0.020em depth0em \kern-0.20em \rm C}
        }}

\newcommand{\R}{\mathbb{R}}

\newcommand{\ba}{\begin{array}}
\newcommand{\ea}{\end{array}}
\newcommand{\be}{\begin{equation}}
\newcommand{\ee}{\end{equation}}
\newcommand{\bea}{\begin{eqnarray}}
\newcommand{\eea}{\end{eqnarray}}
\newcommand{\beaa}{\begin{eqnarray*}}
\newcommand{\eeaa}{\end{eqnarray*}}

\def\d{\delta}

\def\z{\zeta}

\font\tenmath=msbm10 \font\sevenmath=msbm7 \font\fivemath=msbm5
\newfam\mathfam \textfont\mathfam=\tenmath
\scriptfont\mathfam=\sevenmath \scriptscriptfont\mathfam=\fivemath
\def\math{\fam\mathfam}

\def \={{\buildrel {\rm (law)} \over =}}

\def \R{{\math R}}

\def\qed{ \hfill \vrule width.25cm height.25cm depth0cm\smallskip}

\newcommand{\basa}{\begin{assumption}}
\newcommand{\easa}{\end{assumption}}

\newcommand{\bas}{\begin{assum}}
\newcommand{\eas}{\end{assum}}

\newcommand{\ignore}[1]{}
\begin{document}

\renewcommand{\thefootnote}{\fnsymbol{footnote}}

\title{Wiener integrals, Malliavin calculus and covariance  measure structure}
\author{Ida Kruk $^{1}\quad$ Francesco Russo $^{1}\quad$Ciprian A.
Tudor $^{2}\vspace*{0.1in}$\\$^{1}$
 Universit\'{e} de Paris 13, Institut Galil\'ee, Math\'ematiques,\\ 99,
avenue J.B. Cl\'ement, F-93430, Villetaneuse Cedex,
France.\vspace*{0.1in}\\$^{2}$SAMOS/MATISSE, Centre d'Economie de
La Sorbonne,\\ Universit\'e de Panth\'eon-Sorbonne Paris 1,\\90,
rue de Tolbiac, 75634 Paris Cedex 13, France.\vspace*{0.1in}}
\maketitle

\begin{abstract}
We introduce the notion of {\em covariance  measure structure} for
square integrable stochastic processes.
We define  Wiener integral, we  develop a suitable formalism for  
 stochastic calculus of variations and we make Gaussian assumptions only when necessary.
 Our main examples are finite quadratric variation processes with stationary increments and  the bifractional Brownian
motion.
\end{abstract}

\vskip0.5cm



\vfill  \eject

\section{Introduction}
Different approaches have been used to extend the classical It\^o's
stochastic calculus. When the integrator stochastic process does not have the
semimartingale property, then the powerful It\^o's theory cannot be
applied to integrate stochastically. Hence alternative ways have  been
then developed, essentially of two types:
\begin{description}
\item{$\bullet$ } a trajectorial approach, that mainly includes
the rough paths theory (see \cite{QL}) or the stochastic calculus
via regularization (see \cite{RV1993}).
\item{$\bullet$ } a Malliavin calculus (or  stochastic calculus of
variations) approach.
\end{description}
Our main interest consists here in the second approach. 
Suppose that  the integrator is a Gaussian process
$X=(X_{t})_{t\in [0,T]}$. The Malliavin derivation can be naturally
developed on a Gaussian space, see, e.g. \cite{Watanabe},
\cite{Nualart} or \cite{Mall}. A Skorohod (or divergence) integral
can also be  defined as the adjoint of the Malliavin derivative. The
crucial ingredient is the canonical Hilbert space ${\cal{H}}$
(called also, improperly,  by some authors  reproducing kernel Hilbert
space) of the Gaussian process $X$ which is  defined as the closure
of the linear space generated by the indicator functions $\{
1_{[0,t]}, t\in [0,T]\}$ with respect to the scalar product
\begin{equation}
\label{ps} \langle 1_{[0,t]}, 1_{[0,s]}\rangle _{{\cal{H}}}= R(t,s)
\end{equation}
where $R$ denotes the covariance of $X$. Nevertheless, this
calculus remains more or less abstract if the structure of the
elements of the Hilbert  space ${\cal{H}}$ is not known. When we
say abstract, we refer to the fact that, for example, it is
difficult to characterize the processes which integrable with
respect to $X$, to estimate the $L^{p}$-norms of the Skorohod
integrals or to push further this calculus to obtain an
It\^o type formula.

A particular case can be analyzed in a deeper way. We refer here to
the situation when the covariance $R$ can be explicitly written  as
\begin{equation*}
R(t,s)= \int _{0}^{t\wedge s} K(t,u)K(s,u)du,
\end{equation*}
where $K(t,s)$, $0<s<t<T$ is a deterministic kernel satisfying some
regularity conditions. Enlarging, if need, our probability space, we
can express the process $X$ as
\begin{equation}
\label{rep} X_{t}=\int _{0}^{t}K(t,s)dW_{s}
\end{equation}
where $(W_{t})_{t\in [0,T]}$ is a standard Wiener process and the
above integral is understood in the Wiener sense. In this case, more
concrete results can be proved, see \cite{AMN, De, MV}.
In this framework the underlying Wiener process $(W_t)$ is strongly used for 
developing anticipating calculus. In particular \cite{MV} 
puts emphasis on the case $K(t,s) = \varepsilon (t-s)$, when the variance scale 
of the process is as general as possible, including logarithmic scales.

The canonical space ${\cal{H}}$ can be written as
\begin{equation*}
{\cal{H}}= \left( K^{\ast }\right) ^{-1}\left( L^{2}([0,T])\right)
\end{equation*}
where the "transfer operator" $K^{\ast }$ is defined on the set of
elementary functions as
\begin{equation*}
K^{\ast}(\varphi) (s) = K(T,s) \varphi (s)+ \int _{s}^{T} \left(
\varphi (r)-\varphi (s) \right) K(dr,s)
\end{equation*}
and extended (if possible) to
${\cal{H}}$ (or a set of functions contained in ${\cal{H}}$).
Consequently, a stochastic process $u$ will be Skorohod integrable
with respect to $X$ if and only if $K^{\ast } u$ is Skorohod
integrable with respect to $W$ and $\int u\delta X= \int
(K^{\ast}u)\delta W$. Depending on the regularity of $K$ (in
principal the H\"older continuity of $K$ and $\frac{\partial K
}{\partial t}(t,s)$ are of interest) it becomes possible to have
concrete results.

Of course, the most studied case is  the fractional
Brownian motion (fBm), due to the multiple applications of this
process in various area, like telecommunications, hydrology or
economics. Recall that the fBm $(B^{H}_{t})_{t\in [0,T]}$, with
Hurst parameter $H\in (0,1)$ is defined as a centered Gaussian
process starting from zero with covariance function
\begin{equation}
\label{covFBM} R(t,s)= \frac{1}{2}\left( t^{2H}+s^{2H}-\vert
t-s\vert ^{2H}\right), \hskip0.5cm t,s\in [0,T].
\end{equation}
The process $B^{H}$ admits the Wiener integral representation
(\ref{rep}) and the kernel $K$ and the  space ${\cal{H}}$ can be
characterized by the mean of fractional integrals and derivatives,
see \cite{AMN, AN, DU, PiTa, che-nua} among others. As a
consequence, one can prove for any $H$ the following It\^o's formula
\begin{equation*}
f(B^{H}_{t}) = f(0) + \int _{0}^{t} f'(B^{H}_{s})\delta B^{H}_{s} + H\int
_{0}^{t} f''(B^{H}_{s}) s^{2H-1}ds.
\end{equation*}
One can also study the relation between ''pathwise type'' integrals and the
divergence integral, the regularity of the Skorohod integral process
or the It\^o formula for indefinite integrals.

\smallskip

As we mentioned, if the deterministic kernel $K$ in the
representation (\ref{rep}) is not explicitly known, then the
Malliavin calculus with respect to the Gaussian process $X$ remains
in an abstract form; and there are of course many situations when
this kernel is not explicitly known. As main example, we have in
mind the case of the  {\em bifractional Brownian motion (bi-fBm). }
This process, denoted by $B^{H,K}$, is defined as a centered
Gaussian process starting from zero with covariance
\begin{equation}
\label{covbiFBM} R(t,s)= \frac{1}{2^{K}}\left( \left( t^{2H} +
s^{2H}\right) ^{K} -\vert t-s\vert ^{2HK}\right)
\end{equation}
where $H\in (0,1) $ and $K\in (0,1]$. When $K=1$, then we have a
standard fractional Brownian motion.

This process was introduced in \cite{HV} and a ''pathwise type'' approach to stochastic calculus
 was provided in \cite{RT}. An interesting property of
$B^{H,K}$ consists in the expression of its quadratic variation
(defined as usually as limit of Riemann sums, or in the sense of
regularization, see \cite{RV1993}). The following properties hold true.
\begin{description}
\item{$\bullet$ } If $2HK>1$, then the quadratic variation of
$B^{H,K}$ is zero.
\item{$\bullet$ } If $2HK<1$ then the quadratic variation of
$B^{H,K}$ does not exist.
\item{$\bullet$ } If $2HK=1$ then the quadratic variation of
$B^{H,K}$ at time $t$ is equal to $2^{1-K}t$.
\end{description}
The last property is remarkable; indeed, for $HK=\frac{1}{2}$  we
have a Gaussian process which has the same quadratic variation as
the Brownian motion. Moreover, the processes is not a semimartingale
(except for the case $K=1 $ and $H=\frac{1}{2}$), it is
self-similar, has no stationary increments and it is a quasi-helix
in the sense of J.P. Kahane \cite{Kahane}, that is, for all $s\leq
t$,
\begin{equation}
\label{qh}  2^{-K}\vert t-s\vert ^{2HK} \leq E\left| B^{H,K}_{t}-
B^{H,K}_{s}\right| ^{2}\leq 2^{1-K}\vert t-s\vert ^{2HK}.
\end{equation}
We have no information on the form and/or the properties of the
kernel of the bifractional Brownian motion. As a
consequence, a Malliavin calculus was not yet introduced
for this process. On the other side, it is possible to construct a
stochastic calculus of pathwise type, via regularization and one
gets an It\^o formula of the Stratonovich type (see \cite{RT})
\begin{equation*}
f(B^{H,K}_{t}) = f(0)+ \int_{0}^{t} f'(B^{H,K}_{s}) d^{\circ}
B^{H,K}_{s}
\end{equation*}
for any parameters $H\in (0,1)$ and $K\in (0,1]$.

\smallskip

The purpose of this work is to develop a Malliavin calculus with
respect to processes $X$ having a {\em covariance  measure structure} 
in  sense that  the covariance  is the distribution function of
 a (possibly signed) measure on ${\cal{B}}([0,T]^{2})$.
This approach is particularly suitable for processes 
whose  representation form  (\ref{rep}) is not explicitely given.

 We will see  that under this assumption, we can define suitable spaces
on which the construction of the Malliavin derivation/Skorohod
integration is coherent. 

In fact, our initial purpose is more ambitious; we start to
construct a stochastic analysis for general (non-Gaussian) processes
$X$ having a covariance measure $\mu$. We define Wiener integrals
for a large enough class of deterministic functions and we define a
Malliavin derivative and a Skorohod integral with respect to it; we
can also prove certain relations and properties  for these
operators. However,  if one wants 
to produce a consistent theory,   then the Skorohod integral applied to
deterministic integrands should coincide with the Wiener integral.
This property is based on  {\em integration by parts } on  Gaussian spaces which is proved
 in Lemma \ref{l5.4t}. As it can be seen, that proof 
 is completely based on the Gaussian character and
it seems difficult to prove it for general processes.
Consequently, in the sequel, we concentrate  our study on the Gaussian case and
we show various results as  the continuity of the integral
processes, the chaos expansion of local times, the relation between
the ''pathwise'' and the Skorohod integrals and finally we derive the following It\^o formula,
see Corollary \ref{cor8.11c},
for $f \in C^2(\R)$ such that $f''$ is bounded:
$$ f(X_t) = f(X_0) + \int_0^t f'(X_s) \delta X_s + \frac{1}{2} \int_0^t f''(X_s) d\gamma (s), $$
where $\gamma(t) = Var (X_t)$.
 Our main examples include the Gaussian
semimartingales, the fBm with $H \ge  \frac{1}{2}$, the bi-fBm with
$HK\geq \frac{1}{2}$ and processes with stationary increments.
 In the bi-fbm case, when $2HK=1$, we find a very
interesting fact, that is, the bi-fBm with $2HK=1 $ satisfies the
same It\^o formula as the standard Wiener process, that is
\begin{equation*}
f(B^{H,K}_{t})=  f(0) + \int _{0}^{t} f'(B^{H,K}_{s}) \delta
B^{H,K}_{s}+ \frac{1}{2}\int _{0}^{t} f''(B^{H,K}_{s})ds
\end{equation*}
where $\delta $ denotes the Skorohod integral.

We will also like to mention certain aspects that could be the
object of a further study:
\begin{description}
\item{$\bullet$ } the proof of the Tanaka formula involving weighted
local times; for the fBm case, this has been proved in \cite{CNT}
but the proofs necessitates the expression of the kernel $K$.
\item{$\bullet$ } the two-parameter settings, as developed in e.g.
\cite{TV}.
\item{$\bullet$ }the proof of the Girsanov transform and the use of
it to the study of stochastic equations driven by Gaussian noises,
as e.g. in \cite{NO}.
\end{description}

\smallskip

We organized our paper as follows. In Section 2 and 3 we explain the
general context of our study: we define the notion of  covariance measure structure and
we give the basic properties of stochastic processes with this
property. Section 4 contains several examples of processes having
covariance measure $\mu$. Section 5 is consecrated to the
construction of Wiener integrals for a large enough class of
integrands  with respect to (possibly non-Gaussian) process $X$ with
$\mu$. In Section 6, for the same settings, we develop a Malliavin
derivation and a Skorohod integration. Next, we work on a Gaussian
space and our calculus assumes a more intrinsic form; we give concrete
spaces of functions contained in the canonical Hilbert space of $X$
and this allows us to characterize the domain of the divergence
integral, to have Meyer inequalities and other consequences.
Finally, in Section 8 we present the relation ''pathwise''-Skorohod
integrals and we derive an It\^o formula; some particular cases are
discussed in details.

\vskip0.5cm

\section{Preliminaries}
In this paper, a rectangle will be a subset $I$ of $\mathbb{R}^2_+$
of the form
\[
I=]a_1,b_1] \times ]a_2,b_2]
\]
and $T>0$ will be fixed. Given  $F:\mathbb{R}_+ \rightarrow
\mathbb{R}$ we will denote
\[
\Delta_I F=F(b_1,b_2)+F(a_1,a_2)-F(a_1,b_2)-F(b_1,a_2).
\]
Such function will be said to vanish on the axes if
$F(a_1,0)=F(0,a_2)=0$ for every  $a_1, a_2 \in \mathbb{R}_+$.

\vskip0.5cm

Given a continuous function $F:[0,T] \rightarrow \mathbb{R}$ or a process
 $(X_t)_{t \in [0,T]}$, continuous in $L^2(\Omega)$,  will be prolongated  by convention to $\mathbb{R}$ by continuity.

\begin{df} \label{d1.1}
$F:[0,T]^2 \rightarrow \mathbb{R}$ will be said to have a {\bf
bounded planar variation} if
\begin{equation}\label{1.1}
\sup_{ \tau} \sum_{i,j=0}^{n}\left| \Delta_{]t_i,t_{i+1}]\times
 ]t_j,t_{j+1}]}F \right| < \infty.
\end{equation} 
where $\tau  = \{0 = t_0 < \ldots < t_n = 1 \}$ is a subdivsion
of $[0,T]$.
\ A function $F$ will be said to be  {\bf planarly increasing} if for any
rectangle $I \subset  [0,T]^2$  we have $\Delta_I F \geq 0$.
\end{df}

\begin{lem} \label{l1.1}
Let $F:[0,T]^2 \rightarrow \mathbb{R}$ vanishing on the axes having
a bounded planar variation. Then $F=F^+ - F^-$ where $F^+, F^-$ are
planarly increasing and vanishing on the axes.
\end{lem}
{\bf Proof: } It is similar to the result of the one-parameter
result, which states that a bounded variation function  can be
decomposed into the difference of two increasing functions. The
proof of this result is written for instance in \cite{Taylor}
section 9-4. The proof translates into the planar case replacing
$F(b)-F(a)$ with $\Delta_IF$. \qed

\vskip0.5cm

\begin{lem} \label{l1.2}
Let $F:[0,T]^2 \rightarrow \mathbb{R}_+$ be a continuous, planarly
increasing function. Then there is a unique non-atomic, positive,
finite measure $\mu$ on $\mathcal{B}([0,T]^2)$ such that for any $I
\in \mathcal{B}([0,T]^2)$
\[
\mu(I)=\Delta_IF .
\]

\end{lem}
{\bf Proof:  } See Theorem 12.5 of \cite{Bill}. \qed

\vskip0.5cm

\begin{cor}\label{c1.3}
Let $F:[0,T]^2 \rightarrow \mathbb{R}$ vanishing on the axes.
Suppose that $F$ has bounded planar variation. Then, there is a
signed, finite measure $\mu$ on $\mathcal{B}([0,T]^2)$ such that for
any rectangle $I$ of $[0,T]^2$
\[
\Delta_IF=\mu(I).
\]
\end{cor}
{\bf Proof:  } It is a consequence of Lemma \ref{l1.1} and
\ref{l1.2}. \qed

\vskip0.5cm

 We recall now the notion of finite quadratic planar
variation introduced in \cite{RV2000}.

\begin{df} \label{d1.2}
A function $F:[0,T]^2 \rightarrow \mathbb{R}$ has {\bf finite
quadratic planar variation} if
\[
\frac{1}{\varepsilon^2}\int_{[0,T]^2} \left(
\Delta_{]s_1,s_1+\varepsilon] \times ]s_2,s_2+ \varepsilon]} F
\right)^2 ds_1 ds_2
\]
converges. That limit will be called the {\bf planar quadratic
variation of $F$}.
\end{df}

\vskip0.5cm

We introduce now some notions related to stochastic processes. Let
$(\Omega, \mathcal{F},P)$ be a complete probability space. Let
$(Y_t)_{t\in [0,T]}$ with paths in $L^1_{loc}$ and $(X_t)_{t\in
[0,T]}$ be a cadlag  $L^{2}$-continuous process.
Let $t\geq 0$. We denote by
\begin{eqnarray}
I_{\varepsilon}^-(Y,dX,t)&=&\int_0^t Y_s \frac{X_{s+\varepsilon}-X_s}{\varepsilon}ds \nonumber\\
I_{\varepsilon}^+(Y,dX,t)&=&\int_0^t Y_s \frac{X_{s}-X_{s - \varepsilon}}{\varepsilon}ds \nonumber\\
C_{\varepsilon}(X, Y, t) &=& \frac{1}{\varepsilon} \int_0^t 
(X_{s+\varepsilon}-X_s)(Y_{s+\varepsilon}-Y_s)ds. \nonumber
\end{eqnarray}
We set
\[
\int_0^t Y d^-X\ (\textrm{resp.} \  \int_0^t Y d^+X ) 
\]
the limit in probability of
\[
I_{\varepsilon}^-(Y,dX,t)\ (\textrm{resp. } 
I_{\varepsilon}^+(Y,dX,t)).
\]
 $\int_0^tY d^-X$ (resp.  $\int_0^tY d^+X$)
  is called (definite)  {\bf forward} (resp. {\bf backward})
integral of Y with respect to X.
We denote by $[X,Y]_t$ the limit in probability of $C_{\varepsilon}(X, Y,t)$.
$[X,Y]_t$ is called {\bf covariation} of $X$ and $Y$.
If
$X = Y$, $[X,X]$ is called {\bf quadratic variation} of X, also denoted by $[X]$.
\begin{rem} \label{rem2.2} If $I$ is an interval with end-points $a<b$, then
\[
\int_0^T 1_Id^-X = \int_0^T 1_Id^+X = X_b - X_a.
\]

\end{rem}
\vskip0.3cm

Let $(\mathcal{F}_t)_{t \in [0,T]}$ be a filtration satisfying the
usual conditions. We recall, see \cite{RV1993}, that if $X$ is an
$(\mathcal{F}_t)$-semimartingale and $Y$ is a cadlag process (resp.
an $(\mathcal{F}_t)$-semimartingale) then $\int_0^t Y d^-X$ (resp.
$[Y,X]$) is the It\^o integral (resp. the classical covariation).
\par If $X$ is a continuous function and $Y$ is a cadlag function then $\int_0^t Yd^-X$ coincides with the Lebesgue-Stieltjes integral $\int_0^t Y dX$.

\vskip0.5cm

\section{Square integrable processes and covariance measure structure} \label{sec2}

\par In this section we will consider a cadlag zero-mean square integrable process $(X_t)_{t \in [0,T]}$ with covariance
\[
R(s,t) = Cov(X_s, X_t).
\]
For simplicity we suppose that $t \rightarrow X_t$ is continuous in
$L^2(\Omega)$. $R$ defines naturally a finitely additive function
$\mu_R$ (or simply $\mu$) on the algebra  $\mathcal{R}$ of finite
disjoint rectangles included in $[0,T]^{2}$ with values on
$\mathbb{R}$. We set indeed
\[
\mu(I)=\Delta_I R.
\]
A typical example of square integrable processes are Gaussian
processes.

\begin{df} \label{d2.1}
A square integrable process will be said to have a {\bf covariance
measure} if $\mu$ extends to the Borel $\sigma$-algebra
$\mathcal{B}([0,T]^{2})$ to a signed $\sigma$-finite measure.
\end{df}
We recall that $\sigma(I  \mbox{ rectangle, } I\subset [0,T]^{2})
= \mathcal{B}([0,T]^{2} )$.

\vskip0.5cm

\begin{rem} \label{r2.2}
The process $(X_t)_{t\in [0,T]}$ has covariance measure if and
only if $R$ has a bounded planar variation, see Corollary
\ref{c1.3}.\end{rem}

\vskip0.3cm

\begin{df} \label{r2.4}Let us recall   a classical notion introduced in \cite{GR} and
\cite{Follmer}. A process $(X_t)_{t \in [0,T]}$ has {\bf finite
energy} (in the sense of discretization) if
\[
\sup_{\tau } \sum_{i=0}^{n-1}E(X_{t_{i+1}}-X_{t_{i}})^2 < \infty.
\]
Note that if $X$ has a covariance measure then it has finite energy.
Indeed for a given subdivision $t_0<t_1<\ldots<t_n$, we have
\[
\sum_{i=0}^{n-1} E(X_{t_{i+1}}-X_{t_{i}})^2 =
\sum_{i=0}^{n-1}\Delta_{]t_i,t_{i+1}]^2}R\leq \sum_{i,j=0}^{n-1}
\left| \Delta_{]t_i,t_{i+1}] \times ]t_j,t_{j+1}]} R \right| .
\]
\end{df}

\begin{rem} \label{r2.5}
Let $X$ be a process with covariance measure. Then $X$ has a
supplementary property related to the energy. There is a function
${\cal{E}}:[0,T] \rightarrow \mathbb{R}_+$ such that, for each
sequence of subdivisions $(\tau^N)=\left\{ 0= t_0<t_1< \ldots <t_n =
T  \right\}$, whose mesh converges to zero, the quantity
\begin{equation}
    \sum_{i=1}^n E\left( X_{t_{i+1}\wedge t}-  X_{t_{i}\wedge t}\right)^2, \label{2.3}
\end{equation}
converges uniformly in $t$, to ${\cal{E}}$.

 Indeed
\[
\sum_{i=1}^n E\left( X_{t_{i+1}\wedge t}-  X_{t_{i}\wedge
t}\right)^2=\mu(D^N \cap [0,t]^2)
\]
where $ D^N = \bigcup_{i=0}^{n-1}]t_i, t_{i+1}]^2.$  We have $
\bigcap_{N=0}^{\infty}D^N = \left\{(s,s)| s \in [0,T]\right\}. $
From now on we will set
\[
D = \left\{(s,s)| s \in [0,T]\right\} \textrm{ and } D_t = D \cap
[0,t]^2.
\]
Then, for every $t \in [0,T]$, $$ {\cal{E}}(t)=\mu(D_t).$$
\end{rem}

\vskip0.5cm

 We will introduce the  notion of energy in the sense of
regularization, see \cite{RV2000}.
\begin{df} \label{d2.6}
A process $(X_t)_{t \in [0,T]}$ is said to have {\bf finite energy}
if
\[
\lim_{\varepsilon \rightarrow 0}E(C_\varepsilon (X,X,t))
\]
uniformly exists. This limit will be further denoted by
${\cal{E}}(X)(t)$.
\end{df}
\vskip0.3cm

 From now on, if we do not explicit contrary, we will
essentially use the notion of energy in the sense of regularization.

\begin{lem} \label {l2.6}
If $X$ has a covariance measure $\mu$, then it has finite energy.
Moreover
\[
{\cal{E}}(X)(t)=\mu(D_t).
\]
\end{lem}
{\bf Proof:  }It holds that
\begin{eqnarray}
E(C_\varepsilon(X,X,t))&=&\frac{1}{\varepsilon}\int_0^tds E\left( X_{s+\varepsilon}-X_s\right)^2 = \frac{1}{\varepsilon}\int_0^t ds \Delta_{]s,s+\varepsilon]^2}R\nonumber\\
&=& \frac{1}{\varepsilon}\int_0^t ds
\int_{]s,s+\varepsilon]^2}d\mu(y_1,y_2)=
\frac{1}{\varepsilon}\int_{[0,t+\varepsilon ]^2}
d\mu(y_1,y_2)f_{\varepsilon}(y_1,y_2), \nonumber
\end{eqnarray}
where \nopagebreak
\[
f_{\varepsilon}(y_1,y_2)= \left\{
\begin{array}{clrr}
    \frac{1}{\varepsilon}Leb(](y_1-\varepsilon) \vee (y_2 - \varepsilon) \vee 0,y_1 \wedge y_2 ]) & \ \ \mbox{ if }y_1 \in ]y_2-\varepsilon, y_2 + \varepsilon]\\
    0 & \ \ \textrm{otherwise}.
\end{array}
\right.
\]
We observe that
\begin{equation*}
\left| f_\varepsilon (y_1,y_2)\right| \leq \frac{2
\varepsilon}{\varepsilon} = 2 \mbox{ and } f_\varepsilon(y_1, y_2)
\stackrel{\varepsilon \rightarrow 0}{\longrightarrow} \left\{
    \begin{array}{clrr}
    1 & \ \ : y_2 = y_1 \\
    0 & \ \ : y_2 \neq y_1.
\end{array} \right.
 \end{equation*}
So by Lebesgue dominated convergence theorem,
\[
E(C_\varepsilon(X,X,t)) \rightarrow {\cal{E}}(X)(t),
\]
with $ {\cal{E}}(X)(t)=\mu(D_t).  $\qed

\vskip0.5cm

 We recall a result established in \cite{RV2000}, see
Proposition 3.9.
\begin{lem} \label{l2.7}
Let $(X_t)_{t \in [0,T]}$ be a continuous, zero-mean Gaussian
process with finite energy. Then $C_\varepsilon(X,X,t)$ converges in
probability and it is deterministic for every $t \in [0,T]$ if and
only if the planar quadratic variation of $R$ is zero. In that case
$[X,X]$ exists and equals ${\cal{E}}(X)$.
\end{lem}

\vskip0.5cm

This allows to establish the following result.
\begin{prop} \label {p2.8}
Let $(X_t)_{t \in [0,T]}$ be a zero-mean continuous Gaussian process, $X_0=0$,
having a covariance measure $\mu$. Then
\[
[X,X]_t=\mu(D_t).
\]
In particular the quadratic variation is deterministic.
\end{prop}
{\bf Proof: } First, if R has bounded planar variation, then it has
zero planar quadratic variation. Indeed, by Corollary \ref{c1.3} $R$
has a covariance measure $\mu$ and so
\begin{equation}
\frac{1}{\varepsilon} \int_0^1 dt \left(
\Delta_{]t,t+\varepsilon]^2}R \right)^2 \leq \frac{1}{\varepsilon}
\int_0^1 dt \vert \mu \vert \left( ]t,t+\varepsilon]^2 \right)
\Gamma(\varepsilon), \label{2.6}
\end{equation}
where
\[
\Gamma(\varepsilon)= \sup_{|s-t|<\varepsilon} \left|
\Delta_{]s,t]^2}R\right|.
\]
Since $R$ is uniformly continuous, $\Gamma(\varepsilon) \rightarrow
0$. Using the same argument as in the proof of Lemma \ref{l2.6} we
conclude that (\ref{2.6}) converges to zero.

\smallskip

Second, observe that Lemma \ref{l2.6} implies that $X$ has finite
energy. Therefore the result follows from Lemma \ref{l2.7}. \qed

\vskip0.5cm

\section {Some examples of processes with covariance measure} \label{sec3}

\subsection{ $X$ is a Gaussian martingale}

\par It is well known, see \cite{RV2000, Stricker}, that $[X]$ is deterministic. We denote $\lambda(t)=[X]_t$. In this case
\[
R(s_1,s_2) = \lambda(s_1 \wedge s_2)
\]
so that
\[
\mu(B) = \int_B \delta(ds_2-s_1)\lambda(ds_1), \ \ B \in
\mathcal{B}([0,T]^2).
\]
If $X$ is a classical Wiener process, then $\lambda(x)=x$. The support
of $\mu$ is the diagonal $D$, so $\mu$ and the Lebesgue measure are
mutually singular.

\subsection{ $X$ is a Gaussian $(\mathcal{F}_t)$-semimartingale}

\par We recall, see  \cite{Stricker, Emery}, that $X$ is  a semimartingale if and only if it is
a quasimartingale, i.e.
\[
E\left( \sum_{j=0}^{n-1}\left| E\left( X_{t_{j+1}} - X_{t_{j}} | \mathcal{F}_{t_{j}}) \right) \right| \right) \leq K.
\]

\par We remark that if  $X$ is an   $(\mathcal{F}_t)$-martingale or a process such that
  $E\left( \left\|X\right\|_T \right) < \infty$, where    $ \left\|X\right\|_T $
is the total variation, then the
above  condition is easily  verified. According to \cite{Stricker} $\mu$
extends to a measure.

\subsection{ $X$ is a fractional Brownian motion $B^H, H>\frac{1}{2}$}

 We recall that its covariance equals, for every $s_{1}, s_{2} \in
 [0,T]$
\[
R(s_1,s_2) = \frac{1}{2}\left ( s_1^{2H}+s_2^{2H}-\left|
s_2-s_1\right|^{2H}\right).
\]
In that case $\frac{\partial^2R}{\partial s_1 \partial
s_2}=2H(2H-1)\left \vert s_2-s_1\right \vert^{2H-2}$ in the sense of
distributions. Since $R$ vanishes on the axes, we have
\[
R(t_1,t_2) = \int_0^{t_1} ds_1 \int_0^{t_2} ds_2
\frac{\partial^2R}{\partial s_1 \partial s_2}.
\]
The function $R$ has bounded planar variation because
$\frac{\partial^2R}{\partial s_1 \partial s_2}$ is non-negative.
Therefore, for  given $I = ]a_1,a_2] \times]b_1,b_2]$, we have
\[
\left| \Delta _I R \right| = \Delta_I R.
\]
Hence for a subdivision $(t_i)_{i=0}^N$ of $[0,T]$
\[
\sum_{i,j=0}^{N} \left| \Delta_{]t_i,t_{i+1}] \times ]t_j,t_{j+1}]}
R \right| = \sum_{i,j=0}^{N} \Delta_{]t_i,t_{i+1}] \times ]t_j,t_{j+1}]} R = R(T,T).
\]
So the condition (\ref{1.1}) is verified.

\subsection{$X$ is a bifractional Brownian motion with  $H \in (0,1), K \in (0,1]$ and $2HK\geq 1$}

We refer to  \cite{HV}, \cite{RT} for the definition and the basic
properties of this process.  The covariance of the bi-fBm is given
by (\ref{covbiFBM}). We can write its covariance as 
\[
R(s_1,s_2)=R_1(s_1,s_2) + R_2(s_1,s_2),
\]
where
\begin{equation}\label{R1}
R_1(s_1,s_2) = \frac{1}{2^K} \left[ \left(s_1^{2H} + s_2^{2H}\right)^K -
\left( s_1^{2HK} + s_2 ^{2HK}\right) \right] \end{equation} and
\begin{equation}
\label{R2} R_2(s_1,s_2) = \frac{1}{2^K} \left[ -\vert s_2-s_1 \vert^{2HK} + s_1^{2HK}
+ s_2 ^{2HK} \right]. \end{equation} We therefore  have
\[
\frac{\partial^2R_{1}}{\partial s_1 \partial s_2} =
\frac{4H^2K(K-1)}{2^K} \left(s_1^{2H} + s_2^{2H}\right)^{K-2}
s_1^{2H-1}s_2^{2H-1}.
\]
Since $R_1$ is of class $C^2(]0,T]^2)$ and $\frac{\partial^2R_{1}}{\partial s_1 \partial s_2}$
is always negative, $R_1$ is the distribution function of a negative absolutely continuous finite measure, having
 $\frac{\partial^2R_{1}}{\partial s_1 \partial s_2}$ for density.

\smallskip
Concerning the term $R_2$ we suppose $2HK \geq 1$. 
$R_2$ is (up to a constant)
also the covariance function of a fractional Brownian motion of
index $HK$.
\begin{itemize}
    \item If $2HK > 1$ then $R_2$ is the distribution function of an absolutely continuous positive measure with density
$\frac{\partial^2R_{2}}{\partial s_1 \partial s_2} = 2HK (2HK-1)\left| s_1-s_2\right|^{2HK-2}$ which belongs
 of course to $L^1([0,T]^2)$.
 \item If $2HK = 1$, $R_2(s_1,s_2)= \frac{1}{2^{K}}(s_{1}+s_{2}-\vert s_{1}-s_{2}\vert)$.
\end{itemize}

\vskip0.5cm

\subsection{Processes with weakly stationary increments}

A process $(X_t)_{t  \in [0,T]}$ with covariance $R$  is said {\bf
with weakly stationary increments} if for  any $s,t  \in
[0,T]$, $h \geq 0,$ the covariance $R(t+h,s+h)$ does not depend on
h.

\begin{rem} \label{r3.3}
 If $(X_t)_{t  \in [0,T]}$ is a Gaussian process then $(X_t)_{t
\in [0,T]}$ is with weakly stationary increments if and only it
has stationary increments, that is, for  every subdivision $0=t_0
< t_1< \ldots <t_n$ and for every  $h \geq 0$ the law of
$(X_{t_1+h}-X_{t_0+h}, \ldots , X_{t_n+h}-X_{t_{n-1}+h})$ does not
depend on h.
\end{rem}

 We consider a zero-mean continuous in $L^2$ process  $(X_t)_{t  \in [0,T]}$ such that $X_0 = 0$ a.s.
 Let $d(s,t)$ be the associate canonical distance, i.e.
\[
d^2(s,t) = E\left( X_t - X_s \right)^2, \ \ s,t  \in [0,T].
\]
Since $(X_t)_{t  \in [0,T]}$ has stationary increments one can
write
\[
d^2(s,t) = Q(t-s), \ \ \textrm{where } Q(t)=d(0,t).
\]
Therefore the covariance function $R$ can be expressed as
\[
R(s,t)= \frac{1}{2} (Q(s) + Q(t) - Q(s-t)) .
\]
A typical example is provided when $X$ is a fractional Brownian
motion $B^H$. In that case
\[
Q(s) = s^{2H}.
\]

\begin{rem} \label{p3.4}
$X$ has finite energy  if and only if $Q'(0+)$ exists. This follows
immediately from the property
\[
E\left( C_\varepsilon (X,X,t)\right)= \frac{t
Q(\varepsilon)}{\varepsilon}.
\] \qed
\end{rem}

 We can characterize conditions on $Q$ so that $X$ has a covariance measure.
\begin{prop} \label{p3.5}
If $Q''$ is a Radon measure, then $X$
has a covariance measure.
\end{prop}

\begin{rem} \label{r3.6}
Previous assumption is equivalent  to $Q'$ being of  bounded variation. In that case $Q$ is
absolutely continuous.
\end{rem}
{\bf Proof } of the Proposition \ref{p3.5}: Since
\[
R(s_1,s_2) = \frac{1}{2} \left( Q(s_1) + Q(s_2) - Q(s_2 -
s_1)\right)
\]
we have
\[
\frac{\partial^2 R}{\partial s_1 \partial s_2}  =
- \frac{1}{2}\frac{\partial^2}{\partial s_1 \partial s_2} \left(Q(s_2
- s_1)\right) =  \frac{1}{2} Q''(s_2 - s_1)
\]
in the sense of distributions. This means in particular that for
$\varphi, \psi \in \mathcal{D}(\mathbb{R})$ (the space of smooth test
functions with compact support)
\[
 \int_{\mathbb{R}^2} R(s_1,s_2) \varphi'(s_1) \psi'(s_2) ds_1 ds_2 = - \int_{\mathbb{R}} ds_2 \psi(s_2)  \int_{\mathbb{R}} \varphi(s_1) \frac{Q(s_2-ds_1)}{2}.
\] \qed

\vskip0.5cm

\begin{ex} \label{EMFBM}
\par We provide now an example of a process with stationary increments,
investigated for financial applications purposes by \cite{che}.
It is called  {\bf mixed fractional Brownian motion} and it is defined as
 $X = W + B^H$, where $W$
is a Wiener process and $B^H$ is a fractional
Brownian motion with Hurst parameter  $H > \frac{1}{2}$,   independent from  $W$.
\cite{che} proves that $X$ is a semimartingale if and only if 
 $H >  \frac{3}{4}$.

\par $X$ is
a Gaussian process with
\begin{equation}
Q(t) = \left| t \right|+ \vert t \vert^{2H}. \label{3.2}
\end{equation}
Moreover
\[
Q''(dt) = 2 \delta _0 + q'(t)dt,
\]
where $ q(t) = \vert t \vert^{2H-1}2H sign(t).$ 
\end{ex}

 The example above is still very
particular.
\par  Suppose that $Q''$ is a Radon measure. Then, the  function  $Q'$ can be decomposed in the following way
\[
Q' = Q'_{sc}+Q'_d + Q'_{ac},
\]
where $Q'_{sc}$ is continuous and singular, $Q'_{ac}$ is absolutely
continuous and $Q''_d = \sum_{n}\gamma_n \delta_{x_n}$, with $(x_n)$
- sequence of nonnegative numbers and $\gamma_n \in \mathbb{R}$.

\par For instance if $Q$ is as in (\ref{3.2})
then
\begin{equation*}
Q'_{sc} (t) = 0,  \hskip0.3cm Q'_{ac} (t) =  2H\vert t \vert^{2H-1}sign(t),
\hskip0.3cm  Q''_{d}  (t) = 2 \delta _0  .
\end{equation*}

\par A more involved example is the following. Consider Gaussian zero-mean process $(X_t)_{t  \in [0,T]}$ with stationary increments, $X_0 = 0$ such that
\[
Q(t)= \left\{
    \begin{array} {clrr}
    \left| t \right| & \ \ :\left|t\right| \leq \frac{1}{2} \\
    2^{2H-1}\vert t \vert ^{2H}& \ \ : \left|t\right| > \frac{1}{2}
    \end{array}
\right.
\]
In this case it holds that
\begin{equation*}
Q'(t)= \left\{
    \begin{array} {clrr}
    sign(t) & \ \ :\left|t\right| \leq \frac{1}{2} \\
    2^{2H}H \vert t \vert^{2H-1}sign(t)& \ \ : \left|t\right| > \frac{1}{2}
    \end{array}
\right.
\end{equation*} 
and
\begin{equation*}
 Q''_d  =  2 \delta_0 +
(2H-1)\delta_{\frac{1}{2}} + (2H-1)\delta_{-\frac{1}{2}}, 
\end{equation*}
$$
   Q''_{ac} =  
\left \{ 
\begin{array}{ccc}
0 &:& \vert t \vert < \frac{1}{2} \\
2^{2H} H (2H - 1) \vert t\vert ^{2H-2} &:& \vert t \vert > \frac{1}{2}.
\end{array}  
 \right.
$$

\vskip0.3cm

\begin{prop} \label{p}
Let $(X_t)_{t  \in [0,T]}$ be a Gaussian process with stationary
increments such that $X_0 =0$. Suppose that $Q''$ is a measure.
Then
\[
[X]_t = \frac{Q''(\left\{0\right\})t}{2}.
\]
\end{prop}
{\bf Proof: } It follows from the Proposition \ref{p2.8} and the
fact that
\[
\mu(ds_1,ds_2) = ds_1 Q''(ds_2 - s_1).
\] \qed

\vskip0.5cm

\subsection{Non-Gaussian examples}
A wide class of non-Gaussian processes having a covariance 
measure structure can be provided. We will illustrate how to produce such
processes living in the second Wiener chaos. Let us define, for
every $t \in [0,T]$,
\begin{equation*}
Z_{t}= \int _{\mathbb{R}} \left( \int _{0}^{t}  f(u, z_{1})
f(u,z_{2}) du\right) dB_{z_{1}}dB_{z_{2}},
\end{equation*}
where $(B_{z})_{z\in \mathbb{R}}$ is a standard Brownian motion
and $f: \R^2 \rightarrow \R$ is a Borel function such that
\begin{equation} \label{ng}
 \int _{0}^{T} \int _{0}^{T}\left(\int_{\mathbb{R}} f(t,z)f(s,z)
 dz\right) ^{2}dsdt <\infty.
\end{equation}
Now, condition (\ref{ng}) assures that $\frac{\partial
^{2}R}{\partial s \partial t}$ belongs to $L^{1}([0,T]^{2})$.
Clearly, a large class of functions $f$ satisfies (\ref{ng}).

For example, the {\bf Rosenblatt process} (see \cite{taqqu}) is obtained
for $f(t,z)=  (t-z)_{+}^{k-1}$  with $k\in (\frac{1}{4},
\frac{1}{2})$. 
In that case  (\ref{ng}) is satisfied since $k>\frac{1}{4}$.

 The covariance function of the process $Z$ is given by
\begin{eqnarray*}
R(t,s)&=& \int _{\mathbb{R}^{2}} \left( \int _{0}^{t} f(u,z_{1}
)f(u,z_{2}) du\right) \left( \int _{0}^{s} f(v,z_{1} )f(v,z_{2})
dv\right) dz_{1}dz_{2} \\
&=& \int _{0}^{t} \int _{0}^{s} \left(  \int
_{\mathbb{R}^{2}}f(u,z_{1} )f(u,z_{2}) f(v,z_{1}
)f(v,z_{2})dz_{1}dz_{2}  \right) dvdu
\end{eqnarray*}
and thus
\begin{eqnarray*}
\frac{\partial ^{2}R}{\partial s \partial t} &=& \int
_{\mathbb{R}^{2}}f(t,z_{1} )f(t,z_{2}) f(s,z_{1}
)f(s,z_{2})dz_{1}dz_{2}\\
&=& \left( \int _{\mathbb{R}} f(t,z) f(s,z) dz \right) ^{2}.
\end{eqnarray*}
 In
the case of the Rosenblatt process we get $\frac{\partial
^{2}R}{\partial s \partial t} =cst. \vert t-s\vert ^{4K-2} $.

\vskip0.5cm

It is also possible to construct non-continuous processes that
admit a covariance  measure structure. Let us denote by $K$ the
usual kernel of the fractional Brownian motion with Hurst parameter
$H$ (actually, the kernel appearing in the Wiener integral
representation (\ref{rep}) of the fBm) and consider
$(\tilde{N}_{t})_{t\in [0,T]}$ a compensated Poisson process (see
e.g. \cite{KS}). Then $\tilde{N}$ is a martingale and we can define
the integral
\begin{equation*}
Z_{t}= \int _{0}^{t} K(t,s) d\tilde{N}_{s}.
\end{equation*}
The covariance of $Z$ can be written as
\begin{equation*}
R(t,s)= \int _{0}^{t\wedge s}K(t,u)K(s,u) du =\frac{1}{2} \left(
t^{2H} + s^{2H} -\vert t-s\vert ^{2H} \right).
\end{equation*}
Then it is clear that for $H>\frac{1}{2}$ the above process $Z$ has
covariance  measure structure.

\section {The Wiener integral} \label{sec4}

The Wiener integral, for integrators $X$  which are not the classical
Brownian motion, was considered by several authors. Among the most recent references
there are \cite{PiTa} for the case of fractional Brownian motion and
\cite{jolis} when $X$ is a second order process.

We will consider in this paragraph  a zero-mean, square integrable, continuous in $L^2$,
process $(X_t)_{t \in [0,T]}$ such that $X_0 = 0$. We denote by
$R$ its covariance and we will  suppose that $X$ has a covariance
measure denoted by $\mu$ which is not atomic.

We construct here a Wiener integral with respect to  such a process
$X$. Our starting point is the following result, see for instance
\cite{RV2000}:  if $\varphi$ is a bounded variation continuous real
function, it is well known that
\[
\int_0^t \varphi d^-X = \varphi(t)X_t - \int_0^t X_s d \varphi_s,
\ t  \in [0,T].
\]
 Moreover, it holds that
\begin{equation}
\lim_{\varepsilon \rightarrow 0} I^-(\varepsilon, \varphi, dX, t) =
\int_0^t \varphi d^-X \mbox{ in  } L^{2}(\Omega). \label{4.1}
\end{equation}
 We denote by $BV([0,T])$ the space of real functions with bounded variation, defined on $[0,T]$ and by $C^{1}([0,T])$ the set
 of functions on $[0,T]$ of class $C^{1}$. Clearly the above
 relation (\ref{4.1}) holds for $\varphi \in C^{1}([0,T])$.

\smallskip

 By $\mathcal{S}$ we denote  the closed linear subspace of $L^2(\Omega)$ generated by $\int_0^t \varphi d^-X, \varphi \in BV([0,T])$.
 We define $\Phi:BV([0,T]) \rightarrow \mathcal{S}$ by
\[
\Phi (\varphi) = \int_0^T \varphi d^-X.
\]

\smallskip

We introduce the set $L_{\mu }$ as the vector space of Borel
functions $\varphi:[0,T] \rightarrow \mathbb{R}$ such that
\begin{equation}
 \left\| \varphi \right\|_{|\mathcal{H}|}^2 :=  \int_{[0,T]^2} \vert \varphi(u) \vert \vert
\varphi(v)\vert d\left|\mu\right|(u,v) < \infty. \label{spaceLmu}
\end{equation}
We will also use the alternative notation

\begin{equation}
\left\| \varphi \right\|_{|\mathcal{H}|}^2 = \int_{[0,T]^2} \vert \varphi \otimes \varphi \vert d \left| \mu
\right| < \infty. \label{13bis}
\end{equation}

\begin{rem} \label{r5.0}
\begin{description}
\item[a)] A bounded function belongs to $L_{\mu}$, in particular if  $I$ is a real interval, $1_I \in L_{\mu}$.
\item[b)]
If $\phi \in L_\mu$, $1_{[0,t]}\phi \in L_\mu$ for any $t \in [0,T]$.
\end{description}
\end{rem}

\smallskip
For $\varphi,\: \phi \in L_{\mu }$ we set
\begin{equation}
\left\langle \varphi, \phi \right\rangle_{\mathcal{H}} =
\int_{[0,T]^2} \varphi(u) \phi(v) d\mu (u,v).\label{4.2}
\end{equation}

\vskip0.5cm

\begin{lem} \label{l4.1}
Let $\varphi, \; \phi \in BV([0,T])$. Then
\begin{equation}
\left\langle \varphi, \phi \right\rangle_{\mathcal{H}} = E\left(
\int_0^T \varphi d^-X \int_0^T \phi d^-X\right) \label{1}
\end{equation} and
\begin{equation}  \left\|  \varphi \right\|^2_{\mathcal{H}} = E \left (\int_0^T \varphi d^-X \right)^2. \label{esp}
\end{equation}
\end{lem}
{\bf Proof:  } According to (\ref{4.1}), when $\varepsilon
\rightarrow 0$
\begin{equation}
E\left( \int_0^T \varphi(s)
\frac{X_{s+\varepsilon}-X_s}{\varepsilon} ds \int_0^T \phi(u)
\frac{X_{u+\varepsilon}-X_u}{\varepsilon} du \right).\label{2}
\end{equation}
converges to the right member of (\ref{1}). We observe that
(\ref{2}) equals
\begin{eqnarray}
&&\frac{1}{\varepsilon^2}  \int_0^Tds_1 \int_0^Tds_2 \varphi(s_1)
\phi(s_2)E \left( \left(X_{s_1+\varepsilon}-X_{s_1} \right) \left(  X_{s_2+\varepsilon}-X_{s_2}\right)\right) \nonumber \\
&&=  \frac{1}{\varepsilon^2} \int_0^Tds_1 \int_0^Tds_2 \varphi(s_1) \phi(s_2) \Delta _{]s_1,s_s+\varepsilon] \times ]s_2,s_2+\varepsilon]} R\nonumber\\
&&=  \int_{[0,T]^2} d\mu(u_1,u_2) \frac{1}{\varepsilon^2}
\int_{(u_1-\varepsilon)^+}^{u_1} ds_1 \varphi(s_1)
\int_{(u_2-\varepsilon)^+}^{u_2} ds_2 \phi(s_2). \label{*}
\end{eqnarray}
By Lebesgue dominated convergence theorem, when $\varepsilon
\rightarrow 0$, the quantity (\ref{*}) converges to $$\int_{[0,T]^2}
d\mu(u_1,u_2)\varphi(u_1) \phi(u_2)$$ and the lemma is therefore
proved.\qed

\begin{lem} \label{l4.2}
$(\varphi,\phi) \rightarrow \left\langle \varphi,\phi
\right\rangle_\mathcal{H}$ defines a semiscalar product on
$BV([0,T])$.
\end{lem}
{\bf Proof:  } The bilinearity property is obvious. On the other
hand, if $\varphi \in BV([0,T])$,
\begin{equation}
\left\langle \varphi, \varphi \right\rangle_\mathcal{H} = E\left(
\int_0^T \varphi d^-X \right)^2 \geq 0. \label{4.3}
\end{equation}\qed

We denote by $\left\| \cdot \right\|_{\mathcal{H}}$ the associated
seminorm.

\begin{rem} \label{r4.3}
We use the terminology semiscalar product and seminorm  since the
property $\left\langle \varphi, \varphi \right\rangle_\mathcal{H}
\Rightarrow \varphi = 0$ does not necessarily hold. Take for
instance a process
\[
X_t = \left\{
\begin{array} {clrr} 0,& \ \  t \leq t_0 \\
W_{t-t_0},& \ \ t>t_0,
\end{array}
\right.
\]
where $W$ is a classical Wiener process.
\end{rem}

\vskip0.3cm

\begin{rem}
\label{re5.4}
One of the difficulties in the sequel is caused by the fact that $\| \cdot \|_{|\mathcal{H}|}$ does not define a norm. In particular we do not have any triangle inequality.

\end{rem}

\begin{rem}
\label{r5.4bis}
If $\mu$ is a positive measure, then $\| \cdot \|_{|\mathcal{H}|}$ constitutes a true seminorm.
Indeed, if $f \in L_{\mu}$, we have
\begin{eqnarray}
\left\| f \right\|^2_{|\mathcal{H}|} &=& \int_{[0,T]^2} |f(u_1)||f(u_2)|d|\mu|(u_1,u_2)\nonumber \\
 &=& \int_{[0,T]^2} |f|(u_1)|f|(u_2)d\mu(u_1,u_2) = \big\| |f| \big\| ^2_{\mathcal{H}}.\nonumber
\end{eqnarray}
The triangle inequality follows easily.
\end{rem}
\bigskip

\par
In particular if $X$ is a fractional Brownian motion $B$, $H \geq 1/2$, then $\left\| \cdot \right\|_{\mathcal{\vert H \vert}}$ constitutes a norm.

\vskip0.5cm
\par
We introduce the marginal measure $\nu$ associated with $\mu$. We set
\[
\nu(B) = |\mu|([0,T] \times B)
\]
if $B \in \mathcal{B}([0,T])$.

\begin{lem}
\label{le5.4}
If $f \in L_{\mu}$, we have
\[
\left\| f \right\|_{\mathcal{H}} \leq \left\| f \right\|_{|\mathcal{H}|} \leq \left\| f \right\|_{L^2(\nu)},
\]
where $L^2(\nu)$ is the classical Hilbert space of square integrable functions on $[0,T]$ with respect to $\nu$.
\end{lem}
\textbf{Proof: } The first inequality is obvious. Concerning the second one, we operate via Cauchy-Schwartz inequality. Indeed,
\begin{eqnarray}
\left\| f \right\|^2_{|\mathcal{H}|} &=& \int_{[0,T]^2} \vert f (u_1) f (u_2) \vert d |\mu|(u_1,u_2) \nonumber\\
&\leq& \left\{ \int_{[0,T]^2} f^2(u_1) d|\mu| (u_1,u_2) \int_{[0,T]^2} f^2(u_2) d|\mu|(u_1,u_2) \nonumber   \right\}^{\frac{1}{2}} = \int_0^T f^2(u) d\nu(u).
\end{eqnarray} \qed

\par Let $\mathcal{E}$ be the linear subspace of $L_{\mu}$ constituted by the linear combinations $\sum_i a_i 1_{I_i}$, where $I_i$ is a real interval.

\begin{lem} \label{le5.5}
Let $\nu$ be a positive measure on $\mathcal{B}([0,T])$. Then $\mathcal{E}$ (resp. $C^{\infty}([0,T]))$ is dense in $L^2(\nu)$.
\end{lem}
\textbf{Proof: }
\par i) We first prove that we can reduce to Borel bounded functions. Let $f \in L^2(\nu)$. We set $f_n = (f \wedge n) \vee (-n)$. We have $f_n \longrightarrow f $ pointwise (at each point). Consequently the quantity
\[
\int_{[0,T]^2}|f_n-f|^2 d\nu \rightarrow_{n \rightarrow\infty} 0.
\]
by the dominated convergence theorem.

\par ii) We can reduce to simple functions, i.e. linear combination of indicators of Borel sets. Indeed,
 any bounded Borel
function $f$ is the limit of simple functions $f_n$, again pointwise. Moreover the sequence $(f_n)$ can
 be chosen to be bounded by $|f|$.

\par iii) At this point we can choose $f=1_{B}$, where $B$ is a Borel subset of $[0,T]$. By the Radon
property, for every $n$ there is an open subset $O$ of $[0,T]$ with $B \subset  O$, such that $\nu( O \backslash B) < \frac{1}{n}$.
 This shows the existence of a sequence of $f_n = 1_{O_n}$, where $f_n \longrightarrow f$ in $L^2(\nu)$.

\par iv) Since every open set is a union of intervals, if $f = 1_O$, there is a sequence of step
functions $f_n$ converging pointwise, monotonously to $f$.

\par v) The problem is now reduced to take $f=1_I$, where $I$ is a bounded interval.
It is clear that $f$ can be pointwise approximated by a sequence of $C^0$ functions $f_n$ such
that $|f_n| \leq 1 $.

\par vi) Finally $C^0$ functions can be approximated by smooth functions via mollification;
$f_n = \rho_n * f$ and $(\rho_n)$ is a sequence of mollifiers
converging to the Dirac $\delta$-function.

\vskip0.5cm
The part concerning the density of elementary functions is contained in the previous proof.\qed

\vskip0.5cm

We can now establish  an important density proposition.

\begin{prop} \label{pr5.4}
The set $C^\infty ([0,T])$ (resp. $\mathcal{E})$ is dense in $L_\mu$
with respect to  $\left\| \cdot  \right\|_{|\mathcal{H}|}$ and in
particular to the seminorm $\left\| \cdot  \right\|_{\mathcal{H}}$.
\end{prop}

\begin{rem} \label{re5.5}
As observed in Remark \ref{re5.4}, in general $\left\| \cdot  \right\|_{|\mathcal{H}|}$ does not constitute a norm. This is the reason, why we need to operate via Lemma \ref{le5.4}.
\end{rem}
\textbf{Proof } of the Proposition \ref{pr5.4}: Let $f \in L_{\mu}$. We need to find a sequence $(f_n)$ in $C^{\infty}([0,T])$ (resp. $\mathcal{E}$) so that
\[
\left\| f_n - f \right\|_{|\mathcal{H}|} \longrightarrow_{n \rightarrow \infty} 0.
\]
The conclusion follows by Lemma \ref{le5.5} and \ref{le5.4}. \qed

\vskip0.5cm

\begin{cor} \label{c4.6}It holds that
\begin{description}
\item{i) }$\left\langle \cdot, \cdot \right\rangle_\mathcal{H}$ is
a semiscalar product on $L_{\mu }$. \item{ii) } The linear
application
\[
 \Phi: BV([0,T]) \longrightarrow L^2(\Omega)
\]
defined by
\[
\varphi \longrightarrow \int_0^T \varphi d^-X
\]
can be continuously extended to $L_{\mu }$ equipped with the $\left\|  \cdot \right\|_{\mathcal{H}}$-norm.
Moreover we will still have identity (\ref{esp}) for any $\varphi \in L_\mu$.
\end{description}
\end{cor}
{\bf Proof:  } The part i) is a direct consequence of the previous
result. To check ii),  it is only necessary to prove that $\Phi$ is
continuous at zero. This follows from the property (\ref{esp}).
\qed

\vskip0.5cm

\begin{df} We will set $\int_0^T \varphi d X = \Phi(\varphi)$ and it will
be called {\bf the Wiener integral} of $\varphi$ with respect to $X$.
\end{df}
\begin{rem} \label{r4.7}
Consider the relation $\sim $  on $L_{\mu }$ defined by
\[
\varphi \sim \phi \Leftrightarrow \left\|\varphi - \phi
\right\|_{\mathcal{H}} = 0.
\]
Denoting $L^1_\mu$ as the quotient of $L_{\mu }$ through $\sim$ we
obtain a vector space equipped with a true scalar product. However
$L^1_\mu$ is not necessarily complete and so it is not a Hilbert space. For
the simplicity of the notation, we will still denote by $L_{\mu }$
 its quotient with respect to the relation $\sim $.
Two functions $\phi, \varphi$ will be said to be equal in ${L_\mu}$
if $\varphi \sim \phi$.

The fact that $L^1_\mu$  is a metric non complete space, 
does not disturb  the linear extension. The important property
is that $L^2(\Omega)$ is complete.
\end{rem}

\begin{lem}\label{wfb}
Let $h$ be cadlag. Then
\par i) $\int h d^-X=\int h_{-}dX$,
\par ii) $\int h d^+X=\int h dX$,\\
where
\[
h_-(u) = \lim_{s \uparrow  u} h(s).
\]
\end{lem}
{\bf Proof:  } We only prove point i), because the other one behaves similarly. Since $h$ is bounded, we recall by Remark \ref{r5.0} that $h \in L_\mu$. We have
\[\int _{0}^{T}h_{u} \frac{X_{u+\varepsilon }
-X_{u}}{\varepsilon }du = \int _{0}^{T}h^{\varepsilon }dX\]
with
$h^{\varepsilon }_{s}= \frac{1}{\varepsilon } \int _{s-\varepsilon
}^{s} h_{u}du .$ Since
\[
\Vert h^{\varepsilon } - h_-\Vert  _{{\cal{H}}}^2
=\int_{0}^{T} \int _{0}^{T} \left( h^{\varepsilon }(u_1) -h_-(u_1) \right)
\left( h^{\varepsilon }(u_2) -h_-(u_2) \right)d\mu (u_1,u_2)
\]
and $h^{\varepsilon } \rightarrow h_-$ pointwise, the conclusion follows by Lebesgue convergence theorem. \qed

\begin{cor} \label{c5.8} If $h$ is cadlag, then
\[
\int_0^T h d^- X = \int_0^T h dX \ ( = \int_0^T h d^+X).
\]
\end{cor}
{\bf Proof:  } Taking into account Lemma \ref{wfb},
it is enough to show that
\[ \int_0^T (h - h_-)dX = 0.
\]
This follows because
\begin{eqnarray}
\left\| h - h_- \right\|_\mathcal{H}^2 & = & \int\int_{[0,T]^2} (h - h_-)(u_1)(h - h_-)(u_2) d\mu(u_1, u_2)\nonumber\\
& = & \sum_{i,j}(h(a_i) - h(a_{i-}))(h(a_j) - h(a_{j-}))\mu(\{a_i,a_j\}) = 0 \nonumber
\end{eqnarray}
and because $\mu$ is non-atomic. \qed.

\begin{rem} \label{r5.8a}
If $I$ is an interval with end-points $a<b$, then
\[
\int 1_{I} dX = X_b - X_a.
\]
This is a consequence of previous Corollary and Remark \ref{rem2.2}
\end{rem}

\vskip0.5cm

 \begin{ex} \label{EBFBM}
{ \bf The bifractional Brownian motion case: a significant
subspace of $L_{\mu }$. }\end{ex} 

We suppose again that $H K \ge \frac{1}{2}$.
A significant  subspace included in $L_{\mu }$ is the set $L^2([0,T])$
since $X$ is a classical Brownian motion.
If $K=1$ and $H =  \frac{1}{2}$, there is even equality.
 For $K=1$, and $H >  \frac{1}{2}$,  we refer to \cite{AN}.

In the other cases,  we prove that
\begin{equation}
\label{1KH}
 \Vert f\Vert _{\vert {\cal{H}}\vert } ^{2}\leq C(H,K,T) \Vert f\Vert
 _{L^2([0,T])}.
\end{equation}
where $C(H,K)$ is a constant only depending on $H,K$.
It holds that
\begin{eqnarray*}
 \Vert f\Vert _{\vert {\cal{H}}\vert } ^{2}&=& \int _{0}^{T} \int
 _{0}^{T} \vert f(u_1) \vert \vert f(u_2)\vert 
\left \vert\frac{\partial
 ^{2}R}{\partial u_1\partial u_2}(u_1,u_2)  \right \vert du_1du_2 \\
 &\le& C(H,K) \int _{0}^{T} \int
 _{0}^{T} \vert f(u_1) \vert \vert f(u_2)\vert \left( u_1^{2H} +u_2
 ^{2H}\right) ^{K-2}u_1^{2H-1}  u_2^{2H-1}du_1du_2 \\
 &&+ C(H,K)\int _{0}^{T} \int
 _{0}^{T} \vert f(u_1) \vert \vert f(u_2)\vert \vert u_1-u_2\vert ^{2HK-2}
 du_1du_2 \\
 &:=& A + B.
 \end{eqnarray*}

Concerning $B$ we refer to
 \cite{AN}), so we have only to bound the term $A$. 
It gives  
\begin{eqnarray*}
&& \int_{0}^{T} \int_{0}^{T} \vert f(u_{1})\vert \vert f(u_{2})\vert
 \left( u_1^{2H} +u_2
 ^{2H}\right) ^{K-2}u_1^{2H-1}  u_2^{2H-1}
 du_2 du_1 \\
&& \leq C(H,K)\int_{0}^{T} \int_{0}^{T}\vert f(u_{1})\vert \vert
f(u_{2})\vert (u_{1}u_{2})
^{H(K-2) +2H-1} du_{2}du_{1}\\
&&= C(H,K)\left( \int_{0}^{T} \vert f(u)\vert u^{HK-1} du\right)
^{2}\\
 &&\leq C(H,K) \int_{0}^{T} \vert f(u) \vert ^{2} du   \int _{0}^{T}
 u^{2HK-2}du \leq C(H,K,T) \Vert f\Vert ^{2} _{L^{2}([0,T])}.
 \end{eqnarray*}

\vskip0.5cm

Let us  summarize  a few points of our construction. The space
$L_{\mu}$ given by (\ref{spaceLmu}) is, due to Remark \ref{r4.7},
a space with scalar product and it is in general incomplete. The norm
of this space is given by the inner product (\ref{4.2}). We also define
(\ref{13bis}) which is not a norm in general but it becomes a norm
when $\mu$ is a positive measure.

We denote by $\mathcal{H}$ the abstract completion of $L_{\mu }$.

\begin{rem} \label{r5.17bis}
Remark \ref{r5.0} says that $1_{[0,a]}$ belongs to $L_\mu$ for any $a \in [0,T]$. Therefore $\mathcal{H}$ may be seen
 as the closure of $1_{[0,t]}, t \in [0,T]$ with respect to the scalar product
\[
<1_{[0,s]},1_{[0,t]}>_{\mathcal{H}}=R(s,t).
\]
\end{rem}

$\mathcal{H}$ is now a Hilbert space equipped with the scalar
product $\left\langle \cdot,\cdot
\right\rangle_\mathcal{H}$; it coincides with (\ref{1}) when
restricted to $L_{\mu }$. $\mathcal{H}$ is isomorphic to the
self-reproducing kernel space. Generally that space is the space
of $v: [0,T] \rightarrow \mathbb {R}$, $v(t) = E\left(X_t
Z\right)$ with $Z \in L^2(\Omega)$. Therefore, if $Z = \int g dX$,
$g \in L_{\mu }$, we have
\[
v(t) = \int_{[0,t]} \int_{[0,1]} g(s_2)R(ds_1,ds_2).
\]

\begin{prop} \label{p4.8}
Suppose that $X$ is Gaussian. For $f \in L_{\mu }$, we have
\[
E\left(\int_0^T f d X \right)^2 = E\left( \int_0^T f^2(s)d[X]_s +
2 \int_{[0,T]^2}  1_{(s_1 > s_2)}f(s_1)f(s_2)d\mu(s_1,s_2)\right).
\]
\end{prop}
{\bf Proof:  } It is a consequence of the Proposition \ref{p2.8} and Corollary \ref{c4.6}. \qed

\vskip0.5cm

\section{Wiener  analysis for non-Gaussian processes having a covariance measure structure}

 The aim of this section is to
construct some framework of Malliavin calculus for stochastic
integration related to continuous processes $X$, which are $ L^2$-continuous,  with a covariance measure
defined in Section 2 and $X_0 = 0$. We denote by $C_0([0,T])$ the set of continuous functions on
$[0,T]$ vanishing at zero.
In this section we will also suppose that the law of process $X$ on $C_0([0,T])$ has full support, i. e. the probability that
$X$ belongs to any subset of $C_0([0,T])$ is strictly positive.

  We will start with  a general framework.
We will define  the Malliavin derivative with some related properties in this general, not necessarily Gaussian, framework.
A Skorohod integral with respect to $X$ can be
  defined as the adjoint of the derivative operator in some
  suitable spaces. Nevertheless,  Gaussian properties are needed to go into a more practical and
  less abstract situation: for instance if one wants to exhibit concrete
  examples of processes belonging to the domain of the Skorohod
  integral and estimates for the $L^{p} $ norm of the integral. 
  A key point, where the Gaussian nature of the process intervenes is
    Lemma \ref{l5.4t}. We refer also to the comments
  following that lemma.

\smallskip

 We denote by $Cyl$ the set of smooth and cylindrical
random variables of the form
\begin{equation}
F = f\left( \int \phi_1 dX,\ldots, \int \phi_n dX  \right),
\label{5.1}
\end{equation}
where $n \geq 1$, $f \in C_b^{\infty}(\mathbb{R}^n)$ and $\phi_i
\in L_{\mu }$. Here $\int \phi_i dX$ represents the Wiener
integral introduced before Remark \ref{r5.17bis}.

\vskip0.5cm

We denote by $(\mathcal{F}_t)_{t \in [0,T]}$ the canonical
filtration associated with $X$ fulfilling the usual conditions. The
underlying probability space is $(\Omega, \mathcal{F}_T, P)$, where
$P$ is some suitable probability.  For our consideration, it is not
restrictive to suppose that $\Omega = C_0([0,T])$, so that
$X_t(\omega) = \omega(t)$ is the canonical process. We suppose
moreover that the probability measure $P$ has $\Omega$ as support.
According to Section II.3 of \cite{MR}, $\mathcal{F}C^\infty_b$ is
dense in $L^2(\Omega)$, where
\[
\mathcal{F}C^\infty_b = \left\{  f\left(l_1, \ldots, l_m\right),\; m
\in \mathbb{N},\; f\in C^\infty_b(\mathbb{R}^m), \; l_1, \ldots, l_m
\in \Omega^* \right\}.
\]

\smallskip

On the other side, using similar arguments as in \cite{Kuo} one can
prove that for every $l \in \Omega ^*$ there is a sequence of random variables $Z_n
\in \mathcal{S}, Z_n \rightarrow l$ in $L^2(\Omega)$. Thus
$Cyl$ is dense in $L^2(\Omega)$.

\vskip0.5cm

 The derivative operator $D$ applied to $F$ of the form
(\ref{5.1}) gives
\[
DF = \sum_{i=1}^n \partial_i f\left( \int \phi_1 dX, \ldots, \int
\phi_n dX \right) \phi_i.
\]
In particular $DF$ belongs a.s. to $L_{\mu }$ and moreover $E\Vert
DF\Vert _{\vert {\cal{H}}\vert } ^{2} <\infty $.

\smallskip

Recall that the classical Malliavin operator $D$ is an unbounded
linear operator from $L^2(\Omega)$ into $L^2(\Omega; \mathcal{H})$
where $\mathcal{H}$ is the abstract space defined in Section
\ref{sec4}.

\vskip0.5cm

 We define first the set $\left| \mathbb{D}^{1,2} \right|$,
constituted by $F \in L^2(\Omega)$ such that there is a sequence
$(F_n)$ of the form (\ref{5.1}) and there exists  $Z: \Omega
\rightarrow L_{\mu }$ verifying two
conditions:
\begin{description}
    \item[i)]  $\ F_n \longrightarrow F \textrm{ in } L^2(\Omega)$;
    \item[ii)] $ E\Vert DF_{n}-Z\Vert ^{2}_{\vert {\cal{H}}\vert}:=
      E \int _{0}^{T}\int_{0}^{T} \left| D_{u}F_n - Z \right| \otimes \left| D_{v}F_n - Z  \right| d\left| \mu \right|(u,v)
       \stackrel{n \rightarrow \infty}{\longrightarrow} 0.$
\end{description}

\vskip0.5cm

 The set $\mathbb{D}^{1,2}$ will be the vector subspace
of $L^2(\Omega)$ constituted by functions $F$ such that there is
 a sequence $(F_n)$ of the form (\ref{5.1})
\begin{description}
    \item[i)]  $\ F_n \longrightarrow F \textrm{ in } L^2(\Omega)$;
    \item[ii)] $E\left\| DF_n - DF_m \right\| ^{2} _\mathcal{H} \stackrel {n,m \rightarrow \infty}{\longrightarrow} 0$.
\end{description}

\par We will denote by $Z = DF$ the $\mathcal{H}$-valued r.v. such that $\left\|Z -
DF_n\right\|_\mathcal{H} \stackrel {L^2(\Omega)}{\longrightarrow}
0$. If $Z \in L_{\mu }$ a.s. then
\[
\left\|DF\right\|^2_\mathcal{H} = \int_{[0,T]^2}D_{s_1}F D_{s_2}F
d\mu(s_1,s_2).
\]

\smallskip

Note that $\left| \mathbb{D}^{1,2} \right| \subset
\mathbb{D}^{1,2}$ and  $\mathbb{D}^{1,2}$ is a Hilbert space if equipped with the scalar product
\begin{equation}
\left\langle F,G\right\rangle_{1,2} = E(FG) + E \left\langle
DF,DG\right\rangle_{\mathcal{H}}. \label{5.2}
\end{equation}
In general $\left| \mathbb{D}^{1,2} \right|$ is not a linear space equipped with
scalar product since (\ref{13bis}) is not necessarily a norm.

\vskip0.5cm

\begin{rem} \label{r5.3}
$Cyl$ is a vector algebra. Moreover, if $F,\: G \in Cyl$
\begin{equation}
D(F\cdot G) = GDF + FDG . \label{5.5}
\end{equation}
\end{rem}

\vskip0.5cm

We prove some  immediate properties of the Malliavin derivative.
\begin{lem} \label{l5.4}
Let $F \in Cyl$, $G \in \left|\mathbb{D}^{1,2}\right|$. Then $F
\cdot  G \in \left|\mathbb{D}^{1,2}\right|$ and (\ref{5.5}) still
holds.
\end{lem}
{\bf Proof:  } Let $(G_n)$ be a sequence in $Cyl$ such that
\begin{eqnarray}
&&E(G_n - G)^2 \longrightarrow_{n \rightarrow + \infty} 0, \nonumber \\
&&E\left\{ \int_{[0,T]^2}\left|D G_n - DG\right| \otimes \left|D
G_n - DG\right| d\left| \mu \right| \right\} \longrightarrow_{n \rightarrow + \infty}  0.
\label{5.6}
\end{eqnarray}
Since $F \in L^{\infty}(\Omega)$, $FG_n \longrightarrow FG$ in
$L^2(\Omega)$. Remark \ref{r5.3} implies that
\[
D(FG_n)=G_nDF + FDG_n.
\]
So
\begin{equation}
\int_{[0,T]^2}d \left|\mu \right|    \left|G_nDF - GDF\right|
\otimes \left|G_nDF - GDF\right| = \left(G_n - G \right)^2
\int_{[0,T]^2} \left|DF\right| \left|DF\right| d\left|\mu\right|
\label{5.7}
\end{equation}
If $F$ is of type (\ref{5.1}) then
\[
DF = \sum_{i=1}^n Z_i \phi_i,
\]
where $\phi_i \in L_{\mu }$, $Z_i \in L^{\infty}(\Omega)$.
Therefore the expectation of (\ref{5.7}) is bounded by
\begin{equation}
cst. \sum_{i,j=1}^n \int_{[0,T]^2}\left| \phi_i \right| \otimes \left|
\phi_j \right| d\left| \mu \right| E\left( \left( G_n - G \right)^2
Z_iZ_j \right). \label{t}
\end{equation}
When $n$ converges to infinity, (\ref{t}) converges to zero since
$G_n \rightarrow G$ in $L^2(\Omega)$. On the other hand
\begin{eqnarray}
&&\int_{[0,T]^2} d\left| \mu \right| \left( \left| F \left( DG_n - DG \right) \right|
\otimes \left| F \left( DG_n - DG \right) \right| \right) \nonumber \\
&&= \left|F\right|^2 \int_{[0,T]^2} d\left| \mu \right| \left| DG_n
- DG \right| \otimes \left| DG_n - DG \right|. \nonumber
\end{eqnarray}
Since $F \in L^{\infty}(\Omega)$, previous term converges to zero
because of (\ref{5.6}). By additivity the result follows.\qed

\vskip0.5cm

 We denote by $L^2(\Omega; L_{\mu })$ the set of stochastic processes $(u_t)_{t \in [0,T]}$ verifying
\[
E \left( \left\| u  \right\| ^2_{|\cal{H}|} \right) < \infty .
\]
 We can now define the divergence operator (or the  Skorohod
integral) which is an unbounded map defined from $Dom (\delta)
\subset L^2(\Omega; L_{\mu })$ to $L^2(\Omega)$. We say that $u
\in L^2(\Omega; L_\mu)$ belongs to $Dom(\delta )$ if there is a
zero-mean square integrable random variable $Z\in L^{2}(\Omega) $
such that
\begin{equation}
E(FZ) = E\left(\left \langle DF,u\right\rangle_\mathcal{H} \right)
\label{5.3}
\end{equation}
for every $F \in Cyl$. In other words
\begin{equation}
E(FZ) = E\left( \int_{[0,T]^2} D_{s_1}F u(s_2) \mu(ds_1,ds_2)\right)
\mbox{ for every } F \in Cyl.
\end{equation}
Using Riesz theorem we can see that $u \in Dom(\delta)$ if and only
if the map
\[
F \mapsto E\left( \left\langle
DF,u\right\rangle_\mathcal{H} \right)
\]
is continuous linear form with respect to $\left\| \cdot
\right\|_{L^2(\Omega)}$. Since $Cyl$ is dense in $L^2(\Omega)$,
$Z$ is uniquely characterized. We will set
\[
Z = \int_0^T u \delta X.
\]
$Z$ will be called {\bf the Skorohod integral} of $u$ towards $X$.

\begin{df} \label{rem6.3n}
If $u 1_{[0,t]} \in Dom(\delta)$ for any $t \in [0,T]$, we set $\int_0^t u_s 
\delta X_s := \int_0^T u_s 1_{[0,t]} \delta X_s$
\end{df}

\begin{rem} \label{r5.2}
If (\ref{5.3}) holds, then it will be valid by density for every
$F \in  \mathbb{D} ^{1,2} $.
\end{rem}

\vskip0.5cm

 An important preliminary result in the Malliavin calculus is
the following.

\begin{prop} \label{p5.5}
Let $u \in Dom(\delta)$, $F \in \left| \mathbb{D} ^{1,2} \right|$.
Suppose $F \cdot \int_0^T u_s \delta X_s \in L^2(\Omega)$. Then $F u
\in Dom(\delta)$ and
\[
\int_0 ^T F \cdot u_s \delta X_s = F \int_0^T u_s \delta X_s -
\left\langle DF, u \right\rangle_\mathcal{H}.
\]
\end{prop}
{\bf Proof: } We proceed using the duality relation (\ref{5.3}).
Let $F_0 \in Cyl$.
 We need
to show
\begin{equation}
E\left(F_0 \left\{ F \int_0^T u_s \delta X_s - \left\langle
DF,u\right\rangle_{\mathcal{H}} \right\} \right) = E\left(
\left\langle DF_0,Fu \right\rangle_\mathcal{H} \right). \label{5.9}
\end{equation}
Lemma \ref{l5.4} implies that $F_0 F \in \left| \mathbb{D}^{1,2}
\right|$. The left member of
(\ref{5.9}) gives
\begin{eqnarray}
&& E\left( F_0 F \int_0^T u_s \delta X_s \right) - E\left( F_0 \left\langle DF,u\right\rangle_{\mathcal{H}} \right)\nonumber\\
&& = E\left( \left\langle
D\left(F_0F\right),u\right\rangle_{\mathcal{H}} \right) - E\left(
F_0 \left\langle DF,u\right\rangle_{\mathcal{H}} \right) = E ( < D
(F_0F) - F_0 DF,u >_{\mathcal{H}} ). \label{ec1}
\end{eqnarray}
This gives the right member of (\ref{5.9}) because of the Lemma
\ref{l5.4}.  Remark \ref{r5.2} allows to conclude. \qed

\vskip0.5cm

We state a useful Fubini type theorem which allows to interchange
Skorohod and measure theory integrals. When $X$ is a Brownian motion and
the measure theory integral is Lebesgue integral, then the result is
stated in \cite{Nualart}.

\begin{prop} \label{p5.6}
Let $(G, \mathcal{G}, \lambda)$ be a $\sigma$-finite measured space. Let
$u: G \times [0,T] \times \Omega \longrightarrow \mathbb{R}$ be a
measurable random field with the following properties
\begin{description}
    \item[i)] For every $x \in G$, $u(x, \cdot) \in Dom(\delta)$,
    \item[ii)] \[
   E \int_G d\lambda(x_1) \int_G d\lambda(x_2) \int_{[0,T]^2} d \vert \mu \vert \left|u \right|(x_1,\cdot) \otimes \left| u \right| (x_2, \cdot) < \infty,
    \]
    \item[iii)] There is a measurable version in $\Omega \times G$ of the random field
$ \left(\int_0^T u(x,t) \delta X_t \right)_{x \in G}$,

\item[iv)] It holds that
    \[
    \int_G d\lambda(x) E\left(\int_0^T u\left(x,t\right) \delta X_t\right)^2 < \infty.
    \]
\end{description}
Then $ \int_G d\lambda(x) u(x, \cdot) \in Dom(\delta)$  and
\[
\int_0^T \left( \int_G d\lambda(x) u(x, \cdot)  \right) \delta X_t =
\int_G d\lambda(x) \left( \int_0^T u(x,t) \delta X_t \right).
\]
\end{prop}
{\bf Proof: } We need to prove two properties:
\begin{description}
    \item[a)]
    \[
    \int_G d\lambda(x) \left| u \right|(x,\cdot) \in L^2(\Omega; L_{\mu })
    \]
    \item[b)]
    For every $F \in Cyl$ we have
    \[
    E\left(F \int_G d\lambda(x) \int_0^T u(x,t) \delta X_t\right) = E\left( \left\langle DF, \int_G d\lambda(x) u(x,\cdot) \right\rangle_{\mathcal{H}} \right).
    \]
\end{description}
It is clear that without restriction to the generality we can
suppose $\lambda$ to be a finite measure. Concerning  a) we write
\begin{eqnarray} \label{5.10}
E\left( \left | \int_G d\lambda(x) \left|u\right|(x,\cdot) \right |^2_{\mathcal{|H|}} \right )
&& = \int_{[0,T]^2} d \left| \mu \right|(s_1, s_2) \int_G d\lambda(x_1) \left|u\right|(x_1,s_1) \int_G d\lambda(x_2) \left|u\right|(x_2, s_2) \nonumber \\
&& \\
&& = \int_{G \times G} d\lambda(x_1) d\lambda(x_2) \int_{[0,T]} d\left|
\mu \right| (s_1,s_2) \left| u \right|(x_1,s_1) \left| u
\right|(x_2,s_2). \nonumber
\end{eqnarray}
Taking the expectation of (\ref{5.10}), the result a) follows from
ii). For the part  b) let us consider $F \in Cyl$. Classical
Fubini theorem implies
\begin{eqnarray}
&& E\left( F \int_G d\lambda(x)\left( \int_0^T u(x,t) \delta X_t\right) \right) \nonumber \\
&& = \int_G d\lambda(x) E \left( F\int_0^T u(x,t) \delta X_t \right) = \int_G d\lambda(x) E \left( \left\langle DF, u(x, \cdot) \right\rangle_{\mathcal{H}} \right) \nonumber\\
&& = E\left\{ \int_G d\lambda(x) \int_{[0,T]^2} D_{s_1}Fu(x,s_1)D_{s_2}Fu(x,s_2) d\mu(s_1,s_2) \right\} \nonumber\\
&& = E\left( \int_{[0,T]^2} d\mu(s_1,s_2) D_{s_1}F D_{s_2} F \int_G d\lambda(x) u(x,s_1) \int_G d\lambda(x)u(x,s_2) \right) \nonumber\\
&& =  \left\langle DF, \int_G d\lambda(x) u (x, \cdot)
\right\rangle_{\mathcal{H}} \nonumber.
\end{eqnarray}
At this point the proof of the proposition is concluded.\qed

\vskip0.5cm

We denote by $L_{\mu,2}$ the set of $\phi: [0,T]^2 \rightarrow \mathbb{R}$ such that
\begin{itemize}
    \item $\phi(t_1, \cdot) \in L_\mu, \; \forall t_1 \in [0,T]$,
    \item $t_1 \rightarrow \|\phi(t_1, \cdot)\|_{|\mathcal{H}|}\in L_{\mu}$.
\end{itemize}
For $\phi \in L_{\mu, 2}$ we set
\[
\| \phi \|^2_{|\mathcal{H}| \otimes |\mathcal{H}|} = \int_{[0,T]^2} \|\phi(t_1, \cdot)\|_{|\mathcal{H}|} \|\phi(t_2, \cdot)\|_{|\mathcal{H}|} d|\mu|(t_1,t_2).
\]
Similarly to $|\mathbb{D}^{1,2}|$ we will define $|\mathbb{D}^{1,2}(L_{\mu})|$ and even $|\mathbb{D}^{1,p}(L_{\mu})|$, $p \geq 2$.

\vskip0.5cm

\par We first define $Cyl(L_{\mu})$ as the set of smooth cylindrical random elements of the form
\begin{equation} \label{E1}
u_t = \sum_{\ell = 1}^n \psi_\ell (t)G_\ell, \quad t\in [0,T],
\psi_\ell\in L_\mu, G_\ell \in Cyl.
\end{equation}.

On $L_{\mu,2}$ we also define the following inner semiproduct:
\[
<u_1, u_2>_{\mathcal{H} \otimes \mathcal{H}} = \int_{[0,T]^2} <u_1(t_1, \cdot), u_2(t_2, \cdot)>_{\mathcal{H}} d\mu(t_1,t_2).
\]
This inner product naturally induces a seminorm $\| u \|_{\mathcal{H} \otimes \mathcal{H}}$ and we have of course
\[
\| u \|_{|\mathcal{H}| \otimes |\mathcal{H}|} \geq \|u\|_{\mathcal{H} \otimes \mathcal{H}}.
\]

\vskip0.5cm

\par We denote by $|\mathbb{D}^{1,p}(L_{\mu})|$ the vector space of random elements $u: \Omega \rightarrow L_\mu$ such that there is a sequence $(u_n)
 \in Cyl(L_\mu)$ and
\begin{description}
    \item[i)] $ \|u - u_n\|^2_{|\mathcal{H}|} \longrightarrow_{n \rightarrow \infty} 0$ in $L^p(\Omega)$,
    \item[ii)] there is $Z: \Omega \longrightarrow L_{\mu,2}$ with
    $
    \| Du_n - Z \|_{|\mathcal{H}| \otimes |\mathcal{H}|} \longrightarrow 0 \textrm{ in } L^p(\Omega).
    $
\end{description}
We convene here that
\[
Du_n: (t_1,t_2) \longrightarrow D_{t_1}u_n(t_2).
\]

\vskip0.5cm

Note that until now we did not need the Gaussian assumption on $X$. But
this is essential in following result. It says that when the
integrand is deterministic, the Skorohod integral coincide with the
Wiener integral.

\vskip0.5cm

\begin{prop} \label{p5.4b}
Suppose $X$ to be a Gaussian process. Let $h \in L_{\mu }$. Then
\[
\int_0^T h \delta X_s = \int_0^T h dX_s.
\]
\end{prop}
{\bf Proof: } Let $F \in Cyl$. The conclusion follows from the
following Lemma \ref{l5.4t} and density arguments. \qed

\vskip0.5cm

\begin{lem} \label{l5.4t} Let $F \in Cyl$. Then
\begin{equation}
E\left( \left\langle DF,h\right\rangle_\mathcal{H} \right) = E\left(
F \int_0^T h dX  \right).
\end{equation}
\end{lem}
{\bf Proof: } We use the method given  in \cite{Nualart}, Lemma 1.1.
After normalization it is possible to suppose that $\left\| h
\right\|_{\mathcal{H}}=1$. There is $n \geq 1$ such that $F =
\tilde{f}\left( \int h dX, \int \tilde{\phi}_1dX, \dots, \int
\tilde{\phi}_ndX \right)$, $h, \tilde{\phi}_1,\ldots,\tilde{\phi}_1
\in L_{\mu }$, $\tilde{f} \in C^\infty_b (\mathbb{R}^n)$. We set
$\phi_0 = h$ and we proceed by Gram-Schmidt othogonalization. The
first step is given by
\[
Y_1 = \int h dX - \left\langle h, \tilde{\phi}_1
\right\rangle_{\mathcal{H}}\int h dX = \int \phi_1 dX,
\]
where
$
\phi_1 = \frac{h - \left\langle h, \tilde{\phi}_1\right\rangle
h}{\left\| h - \left\langle h, \tilde{\phi}_1 \right\rangle h
\right\|}
$
and so on. Therefore it is possible to find a sequence
$\phi_0, \ldots, \phi_n \in L_{\mu }$  orthonormal with respect to $\left\langle
\cdot, \cdot \right\rangle_{\mathcal{H}}$, such that
\[
F = f\left( \int \phi_0 dX, \ldots, \int \phi_n dX \right),\ f \in
C_b^\infty (\mathbb{R}^{n+1}).
\]
 Let $\rho$ be the density of the standard normal distribution in $\mathbb{R}^{n+1}$,
 i.e. \\$
\rho(x) = (2 \pi)^{-\frac{n+1}{2}}\exp \left(-\frac{1}{2}\sum_{i=0}^n
x_i^2 \right).$  Then we have
\begin{eqnarray}
E\left(\left \langle DF,h\right\rangle_{\mathcal{H}} \right) &&= E\left( \sum_{i = 0}^n \partial_i f\left(\int \phi_0 dX, \dots, \int \phi_n dX \right)  \right) \left\langle \phi_i,h \right\rangle_\mathcal{H} \nonumber\\
&& = E\left( \partial_0 f\left(\int \phi_0 dX, \dots, \int \phi_n dX \right)  \right) = \int_{\mathbb{R}^{n+1}} \partial_0 f(y)\rho(y)dy \nonumber\\
&&  = \int_{\mathbb{R}}f(y)\rho(y)y_0 dy_0 = E\left( F \int h dX
\right) \nonumber
\end{eqnarray}
which completes the proof of the lemma.\qed

\vskip0.5cm

\begin{rem}
It must be pointed out that the Gaussian property of $X$ appears to
be crucial in the proof of Lemma \ref{l5.4t}. Actually we used the
fact that uncorrelated Gaussian random variables are independent
and also the special form of  the derivative of the Gaussian kernel.
As far as we know, there are two possible proofs of this integration
by parts on the Wiener spaces, both using the Gaussian structure:
one (that we used) presented in Nualart \cite{Nualart} and other
given in Bass \cite{bass} using a Bismut's idea and based on the
Fr\'echet form of the Malliavin derivative.
\end{rem}

\section{The case of  Gaussian processes with a
covariance measure structure}

 Let $X = (X_t)_{t \in [0,T]}$ be a zero-mean Gaussian process
such that $X_0 = 0$ a.s.  that is continuous. A classical result of
\cite{Fernique} (see Th. 1.3.2. and Th. 1.3.3) says that
\begin{eqnarray} \label{EsupG}
 \sup_{t \in [0,T]}  \vert X \vert \in L^2.
\end{eqnarray}
This implies in particular that $X$ is
 $L^2$-continuous.
We suppose also as in previous section that
the law of $X$ in $C_0 ([0,T])$ has full support.

We suppose moreover that {\em it has
covariance $R$ with covariance measure $\mu$.}
 Since $X$ is Gaussian, according to the Section \ref{sec4}, the canonical Hilbert  space ${\cal H}$  of
 $X$ (called reproducing kernel Hilbert space by some authors)
provides an abstract Wiener space and an abstract
  structure of Malliavin calculus was developed, see for instance \cite{Shig, NualartStF,
 Watanabe}.

Recently, several papers were written in relation to fractional Brownian motion and
Volterra  processes of the type
$X_t = \int_0^t G(t,s) dW_s$, where $G$ is a deterministic kernel, see for instance \cite{AMN, De}.
In this work we remain close to the intrinsic approach based on the covariance as in
 \cite{Shig, NualartStF,  Watanabe}.
However their approach is based on a version of self-reproducing kernel space ${\cal H}$ which
is  abstract.
Our construction  focuses on the linear subspace $L_\mu$ of ${\cal H}$
which is constituted by functions.

 \subsection{Properties of Malliavin derivative and divergence operator}

We introduce  some elements of the Malliavin calculus with respect to
$X$. Remark \ref{r5.17bis} says  that the abstract Hilbert space ${\cal H}$ introduced in Section 5
is the topological  linear space generated by the indicator
functions $\{ 1_{[0,t]}, t\in [0,T]\}$ with respect to the scalar
product
\begin{equation*}
\langle 1_{[0,t]}, 1_{[0,s]}\rangle _{\cal{H}}= R(t,s).
\end{equation*}
In general, the elements of ${\cal{H}}$ may not be functions
but distributions. This is actually the case of the fractional
Brownian motion with $H>\frac{1}{2}$, see Pipiras and Taqqu
\cite{PiTa}. Therefore it is more convenient  to work with
one subspace of ${\cal{H}}$ that contains only functions, for instance
  $L_{\mu }$.

\vskip0.5cm

We establish here some peculiar and useful properties of Skorohod integral.

\begin{prop} \label{pro2}
Let $u \in Cyl(L_{\mu})$. Then $u \in Dom(\delta)$ and $\int_0^T u \delta X \in L^p(\Omega) $ for every $1 \leq p < \infty$.
\end{prop}
{\bf Proof: }
Let  $u = G \psi$, $\psi \in L_{\mu}$, $G \in Cyl$. 
Proposition \ref{p5.4b} says that $\psi \in Dom(\delta)$. Applying Proposition \ref{p5.5}
with $F = G$ and $u = \psi$, we get that $\psi G$ belongs to $Dom(\delta)$ and
$$ \int_0^T u \delta X = G \int_0^T \psi \delta X_s - \int_{[0,T]^2} \psi(t_1) D_{t_2} G d\mu(t_1,t_2). $$

If  $G = g(Y_1, \ldots, Y_n)$, where  $Y_i = \int \phi_i dX, \ 1 \leq i \leq n$,
then
\begin{equation} \label{E11}
 \int_0^T u \delta X = - \sum_{j= 1}^n <\phi_j,\psi>_{\mathcal{H}}  \partial_j g(Y_1, \ldots, Y_n)+ g(Y_1, \ldots, Y_n) \int_0^T \psi dX.
\end{equation}
The right member belongs obviously to each $L^p$ since $Y_j$ is a Gaussian random variable and
$g, \partial_j g$ are bounded.
The final result for $u \in Cyl(L_\mu)$ follows by linearity.
\qed

\begin{rem} \label{Rem2}
(\ref{E11}) provides an explicit expression of $\int_0^T u \delta X$.
\end{rem}

\vskip0.5cm

We discuss now the commutativity of the derivative and Skorohod integral. First we observe that if $F \in Cyl$, $(D_t F)\in Dom(\delta)$.
 Moreover, if $u \in Cyl (L_{\mu})$, $(D_{t_1}u(t_2))$ belongs to $|\mathbb{D}^{1,2}(L_{\mu,2})|$. 
Similarly to (1.46) Ch. 1 of \cite{Nualart}, we have the following property.
\begin{prop} \label{P3}
Let $u \in Cyl(L_{\mu})$. Then
\[
\int_0^T u \delta X \in |\mathbb{D}^{1,2}|
\]
and we have 
 for every $t$
\begin{equation} \label{E11*}
D_t \left( \int_0^T u \delta X \right) = u_t + \int_0^T(D_tu_s)\delta X_s.
\end{equation}
\end{prop}
{\bf Proof: } It is enough to write the proof for $u = \psi G$, where $G \in Cyl$ of the type
\[
G = g(Y_1, \ldots, Y_n),
\]
$Y_i = \int \phi_i dX$, $1 \leq i \leq n$. According to (\ref{E11}) in the proof of Proposition \ref{pro2}, the left member of (\ref{E11*}) gives
\begin{eqnarray} \label{E15}
- \sum_{j= 1}^n &&<\phi_j,\psi>_{\mathcal{H}}  D_t (\partial_j g(Y_1, \ldots, Y_n))+ D_t G \int_0^T \psi dX + G \psi (t) \nonumber \\
&=&  - \sum_{j= 1}^n <\phi_j,\psi>_{\mathcal{H}}  \sum_{l=1}^n \partial_{il}^2  g(Y_1, \ldots, Y_n)\phi_l(t)   \\
&& \ \ \ \ \ \ + \sum_{j=1}^{n}  \partial_j g(Y_1, \ldots, Y_n) \int_0^T \psi dX \phi_j(t) + g(Y_1, \ldots, Y_n) \psi (t) .\nonumber
\end{eqnarray}
On the other hand
\[
D_tu(s) = \psi(s) \sum_{j = 1}^{n} \partial_j g (Y_1, \ldots, Y_n) \phi_j(t).
\]
Applying again (\ref{E11}), through linearity, we obtain for $t \in [0,T]$,
\begin{eqnarray}
\int_0^T D_t u \delta X &=& \sum_{j=1}^{n} \phi_j (t) \int_0^T \psi \partial_j g(Y_1, \ldots, Y_n) \delta X \nonumber\\
&=& \sum_{j=1}^{n} \phi_j (t) \left[- \sum_{l = 1}^n <\phi_l,\psi>_{\mathcal{H}} \partial^2_{lj}g(Y_1, \ldots, Y_n) + \partial_j g(Y_1, \ldots, Y_n) \int_0^T \psi dX
 \right]. \nonumber
\end{eqnarray}
Coming back to (\ref{E15}) we get
\[
D_t\left( \int_0^T \psi G \delta X \right) = \int_0^T D_t(\psi G) \delta X + \psi(t)G.
\]\qed

\vskip0.5cm
We can now evaluate the $L^{2}(\Omega)$ norm of the Skorohod integral.

\begin{prop} \label{P4}
Let $u \in |\mathbb{D}^{1,2}(L_{\mu})|$. Then $u \in Dom(\delta)$, $\int_0^T u \delta X \in L^2(\Omega)$ and
\begin{equation} \label{E13}
E\left( \int_0^T u \delta X \right)^2 = E(\left\| u \right\|_{\mathcal{H}}^2) + E\left( \int_{[0,T]^2} d \mu (t_1,t_2)
\int_{[0,T]^2}d\mu(s_1,s_2) D_{s_1}u_{t_1} D_{t_2}u_{s_2}  \right).
\end{equation}
Moreover
\begin{equation} \label{E14}
E\left(\int_0^T u \delta X\right)^2 \leq E \left(  \left\| u \right\|^2_{|\mathcal{H}|} +
 \int_{[0,T]^2}d \vert \mu \vert (t_1,t_2) \left\| D_{\cdot}u_{t_1} \right\|^2_{|\mathcal{H}|}  \right).
\end{equation}
\end{prop}

\begin{rem} \label{Re4}
\begin{description}
\item[i)]
Let $u,v \in |\mathbb{D}^{1,2}(L_{\mu})|$. Polarization identity implies
\begin{eqnarray} \label{E16}
E \left(  \int_0^T u \delta X \int_0^T v \delta X \right) &=& E\left(  <u,v>_\mathcal{H} \right)\nonumber\\
\\ && \ \ \ \ \ + E \left(  \int_{[0,T]^2} d\mu(t_1,t_2) \int_{[0,T]^2}d\mu(s_1,s_2) D_{s_1}u_{t_1} D_{t_2}v_{s_2} \right).\nonumber
\end{eqnarray}
\item[ii)]
If $u \in |  \mathbb {D}^{1,2}(L_\mu)|$, then $u 1_{[0,t]} \in|  \mathbb{D}^{1,2}(L_\mu)| $ for any $t \in [0,T]$ 
and consequently to $Dom(\delta).$
\end{description}
\end{rem}
{\bf Proof }(of Proposition \ref{P4}):
Let $u \in Cyl(L_{\mu})$. By the Proposition \ref{P3}, since $\int_0^T u \delta X \in \mathbb{D}^{1,2}$ we get
\begin{eqnarray}
&E&\left( \int_0^T u \delta X \right)^2 = E\left(  \left\langle u, D\int_0^T u \delta X \right\rangle_{\mathcal{H}} \right) \nonumber\\
&& \ \ \  = E\left( \int_{[0,T]^2} u_{t_1} D_{t_2}\left(\int_0^T u \delta X\right) d\mu(t_1,t_2)  \right) \nonumber\\
&& \ \ \ \ \ \ \  = E\left( \int_{[0,T]^2} u_{t_1} \left( u_{t_2} + \int_0^T D_{t_2} u_s \delta X_s\right) d\mu(t_1,t_2)  \right) \nonumber\\
&& \ \ \ \ \ \ \ \ \ =E\left( \left\| u \right\|^2_{\mathcal{H}} \right) + \int_{[0,T]^2} d\mu(t_1,t_2) E\left( u_{t_1} \int_0^T D_{t_2} u_s \delta X_s\right). \nonumber
\end{eqnarray}
Using again the duality relation, we get
\[
E\left( \left\| u \right\|^2_{\mathcal{H}} + \int_0^T d\mu(t_1,t_2)<D_{\cdot}u_{t_1},D_{t_2}u_{\cdot}>_{\mathcal{H}}\right),
\]
which constitutes formula (\ref{E13}).

Moreover, using Cauchy-Schwarz, we obtain
\begin{equation} \label{E16a}
E\left( \int_0^T u \delta X \right)^2 \leq E\left( \left\| u \right\|^2_{\mathcal{H}} \right) +
E\left( \int_{[0,T]^2} d|\mu|(t_1,t_2) \left\|D_{t_1}u_\cdot \right\|_{\mathcal{H}} \left\|D_{\cdot}u_{t_2}\right\|_{\mathcal{H}} \right).
\end{equation}
Since
\[
\int_{[0,T]^2} d|\mu|(t_1,t_2) \left\|Du_{t_2}\right\|^2_{\mathcal{H}} = \int_{[0,T]^2} d|\mu|(t_1,t_2) \left\|D_{t_1}u\right\|^2_{\mathcal{H}}
\]
(\ref{E16a})  is equal or smaller than
\[
E\left( \int_{[0,T]^2} d|\mu|(t_1,t_2) \left\|Du_{t_2}\right\|^2_{\mathcal{H}} \right)
\]
and this shows (\ref{E14}).
\par Using the fact that $Cyl(L_{\mu})$ is dense in $|\mathbb{D}^{1,2}(L_\mu)|$ we obtain the result. \qed
\vskip0.5cm

\subsection{Continuity of the Skorohod integral process}

\par  It is possible to connect previous objects with the classical Wiener analysis on an abstract Wiener space, related to Hilbert spaces $\mathcal{H}$,
 see \cite{NualartStF, Shig}.
\par In the classical theory the Malliavin gradient (derivative) $\nabla$ and the divergence operator $\delta$ are well defined with its domain. For instance $\delta: \mathbb{D}^{1,2}(\mathcal{H})\rightarrow L^2(\Omega)$ is continuous and $\mathbb{D}^{1,2}(\mathcal{H})$ is contained in the classical domain. However as we said the realizations of $u \in \mathbb{D}^{1,2}(\mathcal{H})$ may not be functions.
\par If $u \in |\mathbb{D}^{1,2}(L_{\mu})|$, it belongs to $\mathbb{D}^{1,2}(\mathcal{H})$ and its norm is given by
\[
\left\| u \right\|^2_{1,2} = E\left(  \left\|u\right\|^2_\mathcal{H} + \int_{[0,T]^2}d\mu(s_1,s_2) \left\|D_\cdot u_{s_1}\right\|_{\mathcal{H}} \left\|D_\cdot u_{s_2}\right\|_{\mathcal{H}}\right).
\]

Classically $\nabla u$ is an element of $L^2(\Omega, \mathcal{H} \otimes \mathcal{H})$, where $\mathcal{H} \otimes \mathcal{H}$ is the Hilbert space of bilinear continuous functionals $\mathcal{H} \otimes \mathcal{H} \rightarrow \mathbb{R}$ equipped with the Hilbert-Schmidt norm.
\par Given $u \in |\mathbb{D}^{1,2}(L_\mu)| \subset \mathbb{D}^{1,2}(\mathcal{H})$, we have $Du \in L^2({\Omega; L_{\mu,2}})$. The associated gradient $\nabla u$ 
is given by
\[
(h,k) \mapsto \int_{[0,T]^2}D_{s_1}u_{t_1}h(s_2)k(t_2)d\mu(s_1,s_2)d\mu(t_1,t_2),
\]
where $h,k \in L_\mu$. Its Hilbert-Schmidt norm coincides with
\[
\int_{[0,T]^2}<Du_{s_1},Du_{s_2}>_{\mathcal{H}} d\mu(s_1,s_2).
\]

\begin{rem} \label{R72}
If $u \in |\mathbb{D}^{1,2}(L_\mu)|$
\[
\int_0^T u_s dX_s = \delta(u).
\]
\end{rem}

\begin{rem} \label{R73}
The standard Sobolev-Wiener space $\mathbb{D}^{1,p}(\mathcal{H})$, $p \geq 2$ is included in the classical domain of $\delta$ and the  Meyer's inequality holds:
\begin{equation} \label{meyer}
E\vert \delta (u)\vert ^{p} \leq C(p) E \left( \Vert u \Vert^{p} _{ {\mathcal{H}} }+ \Vert \nabla u \Vert ^{p}_{   {\mathcal{H}}  \otimes  {\mathcal{H}}}
\right) .
\end{equation}
This implies that if $u \in |\mathbb{D}^{1,2}(L_\mu)|$
\begin{equation} \label{meyer1}
E\left| \int_0^T u \delta X \right|^p \leq C(p)E \left( \left\| u \right\|^p_{\mathcal{H}}  + \left\{\int_{[0,T]^2} <Du_{s_1},Du_{s_2}>_{\mathcal{H}}d\mu(s_1,s_2) \right\} ^{\frac{p}{2}}\right).
\end{equation}
Consequently this gives
\begin{equation} \label{meyer2}
E \left( \left \vert \int_0^T u \delta X \right \vert^p \right) \leq  C(p)E
\left\{ \int_{[0,T]^2} \left|<Du_{s_1},Du_{s_2}>_{\mathcal{H}} \right| d|\mu|(s_1,s_2) \right\}^{\frac{p}{2}}. 
\end{equation}
The last inequality can be shown in a similar way as in the case of  Brownian motion.
One applies Proposition 3.2.1 p. 158 in \cite{Nualart}
and then one argument in the proof of  Proposition  3.2.2 again in \cite{Nualart}.

\end{rem}

The Meyer inequalities are very useful in order to prove the
continuity of the trajectories for Skorohod integral processes. We
illustrate this in the next proposition.

\begin{prop}
Assume that the covariance measure of the process $X$ satisfies
\begin{equation}
\label{cont1} \left[ \mu \left( (s,t]\times (s,t]\right) \right]
^{p-1} \leq \vert t-s\vert ^{1+\beta }, \hskip0.5cm \beta >0
\end{equation}
for some $p>1$ and consider a process $u\in |\mathbb{D}^{1,2p}(L_{\mu
})|$ such that
\begin{equation}
\label{cont2} \int _{0}^{T} \int_{0}^{T} \left(  \int _{0}^{T}
\int_{0}^{T} \vert D_{a}u_{r} \vert \vert D_{b}u_{\theta }\vert
d\vert \mu \vert (a,b) \right) ^{p} d\vert \mu \vert (r,\theta )
<\infty.
\end{equation}
Then the Skorohod integral process $\left( Z_{t}=\int_{0}^{t}
u_{s}\delta X_{s}\right) _{t\in [0,T]} $ admits a continuous
version.
\end{prop}
{\bf Proof: } We can assume that the process $u$ is centered
because the process $\int_{0}^{t}E(u_{s})\delta X_{s}$ always
admits a continuous version  under our hypothesis. By
(\ref{meyer2}), (\ref{cont1}) and (\ref{cont2}) we have
\begin{eqnarray*}
E\left| Z_{t}-Z_{s} \right| ^{2 p} &\leq & c(p) E \left(
\int_{0}^{T}\int_{0}^{T}\int_{0}^{T}\int_{0}^{T} \left| D_{a}
u_{r}1_{(s,t]}(r) \right| \left| D_{b}u_{\theta }1_{(s,t]}(\theta
)\right| d\vert \mu \vert (a,b) d\vert \mu \vert (r,\theta ) \right) ^{p}\\
&\leq & c(p)\left[ \mu \left( (s,t]\times (s,t]\right) \right]
^{p-1}\int _{0}^{T} \int_{0}^{T} \left(  \int _{0}^{T} \int_{0}^{T}
\vert D_{a}u_{r} \vert \vert D_{b}u_{\theta }\vert d\vert \mu \vert
(a,b) \right) ^{p} d\vert \mu \vert (r,\theta )\\
&\leq & c(p,T) \vert t-s\vert ^{1+\beta }.
\end{eqnarray*}
 \qed

\vskip0.5cm

\begin{rem}
In the fBm case we have that
$$\mu \left( (s,t]\times (s,t]\right) = \vert t-s\vert ^{2H}$$
and (\ref{cont1}) holds with $pH>1$.  In the bifractional case, it
follows from \cite{HV} that
$$\left| \mu \left( (s,t]\times (s,t]\right) \right| \leq 2^{1-K}
\vert t-s\vert ^{2HK}$$ and therefore (\ref{cont1}) holds if
$pHK>1$.
\end{rem}

\subsection{On local times }

We will make in this paragraph  a few observations on the chaotic
expansion of the local time  of a Gaussian process $X$ having a
covariance  measure structure. Our analysis is basic and we will
only aim to anticipate  a possible further study. We illustrate the
fact that the covariance measure appears to play an important
role for the existence and the regularity of the local time.

Let  us use the standard way to introduce the local time $L(t,x)$ of
the process $X$; that is, for every $t\in [0,T]$ and $x\in
\mathbb{R}$, $L(t,x)$ is defined as the density of the occupation
measure
\begin{equation*}
\mu _{A}=  \int _{0}^{t} 1_{A}(X_{s}) ds, \hskip0.5cm A \in
{\cal{B}}(\mathbb{R}).
\end{equation*}
It can be formally written as
$$L(t,x) =\int_{0}^{t} \delta _{x}(X_{s})ds,$$
where $\delta $ denotes the Dirac delta function and $\delta
_{x}(X_{s})$ represents a distribution in the Watanabe sense,
see \cite{Watanabe}.

\smallskip

Since $X$  is a Gaussian process, it is possible to construct related  multiple
Wiener-It\^o integrals. We refer to \cite{Nualart} or
\cite{Major} for the elements of this construction.

There is a standard method to compute the Wiener-It\^o chaos
expansion of $L(t,x)$. It consists in approaching the Dirac function by
 mean-zero Gaussian kernels $p_{\varepsilon }$ of variance $\varepsilon$  and
to take the limit in $L^{2}(\Omega)$ as $\varepsilon \to 0$. We get
(see e.g. \cite{Edd})
\begin{equation}
\label{chaos} L(t,x) =\sum _{n=0}^{\infty }\int _{0}^{t} \frac{
p_{R(s,s) }(x)}{ R(s,s) ^{\frac{n}{2}}}H_{n}\left(
\frac{x}{\sqrt{R(s,s)}} \right) I_{n}\left( 1_{[0,s ]}^{\otimes
n}\right) ds
\end{equation}
for all $t\in [0,T]$, $x\in \mathbb{R}$ where $I_{n}$ denotes the multiple Wiener integral of order $n$ with respect to $X$ and
$H_{m}$ represents the Hermite polynomial of degree $m$.
For recent references on the subject one can see \cite{NualartStF, LN}.
 Note that the integral $I_{1}$ is nothing else that the Wiener
integral discussed in Section 5. One can compute the second moment of $L(t,x)$ by using the isometry of multiple stochastic
integrals
\begin{equation*}
E I_{n}\left( 1_{[0,s ]}^{\otimes n}\right)I_{m}\left( 1_{[0,t
]}^{\otimes m}\right) = \left \{
\begin{array}{ccc}
m!R(s,t)^{m} &\mbox{ if }&  m=n \\
0   &\mbox{ if }&  m \neq n
\end{array}
\right.
\end{equation*}
 Using standard bounds as in \cite{Edd},
it follows that the $L^{2}$ norm of (\ref{chaos}) is finite if
\begin{equation*}
\sum_{n\geq 1} n^{-\frac{1}{2}} \int _{0}^{t} \int_{0}^{t} \frac{
\left| \mu \left( [0,u]\times [0,v]\right) \right| ^{n} } {\left(
\left| \mu \left( [0,u]\times [0,u]\right) \right| \left| \mu \left(
[0,v]\times [0,v]\right) \right| \right)
^{\frac{n+1}{2}}}dvdu<\infty.
\end{equation*}
It can be seen that the existence of the local time $L(t,x)$ as
random variable in $L^{2}(\Omega)$ is closely related to the
properties of the covariance measure $\mu$. A possible condition to
ensure the existence of $L$ could be
\begin{equation*}
\int _{0}^{t} \int_{0}^{t} \frac{ \left| \mu \left( [0,u]\times [0,v]\right) \right| ^{n} } {\left( \left| \mu \left(
[0,u]\times [0,u]\right) \right| \left| \mu \left( [0,v]\times [0,v]\right) \right| \right) ^{\frac{n+1}{2}}}dvdu < {\rm cst.}
n^{-\beta }
\end{equation*}
with $\beta > \frac{1}{2}$. Of course, this remains rather abstract
and it is interesting to be checked in concrete cases. We refer to
\cite{NV} for the Brownian case, to \cite{Edd} for the fractional
Bownian case and to \cite{RT} for the bifractional case.

\smallskip

We also mention that the properties of the covariance measure of
Gaussian processes are actually crucial to study sample path
regularity of local times like level sets, Hausdorff dimension etc.
in the context of the existence of {\em local non-determinism. } We
refer e.g. \cite{xiao} for a complete study of path
properties of Gaussian random fields and to \cite{xiao-tud}
for the case of bifractional Brownian motion.

\vskip0.5cm

\section{It\^o formula in the Gaussian case and related topics }

 The next step will consist in expressing the relation
between Skorohod integral and integrals obtained via regularization.
The first result is illustrative. It does not enter into specificity
of the examples.

\vskip0.3cm

\begin{thm} \label{t5.7}
Let $(Y_t)_{t \in [0,T]}$ be a cadlag process. We take into
account the following hypotheses
\begin{description}
    \item[a)] $\sup_{t \leq T} \left| Y_t \right|$ is square integrable.
    \item[b)] $Y \in |\mathbb{D}^{1,2}(L_{\mu})|$. Moreover $DY$ verifies
         \begin{equation}
    \left|D_{t_1} Y_{t_2}\right| \leq Z_2,\ \ \forall (t_1,t_2)\in [0,T]^{2}\ \  |\mu| \;a.e. \label{5.11}
    \end{equation}
    where $Z_2$ is a square integrable random variable.
    \item[c)] For $|\mu|$ almost all $(t_1,t_2) \in [0,T]^2$
    \begin{equation}
    \lim_{\varepsilon \rightarrow 0} \frac{1}{\varepsilon} \int_{t_2 - \varepsilon}^{t_2} D_{t_1}Y_s ds \label{5.12}
    \end{equation}
exists a.s. This quantity will be denoted $(D_{t_1}Y_{t_2-})$. Moreover for each $s$, the set of $t$ such that $D_tY_{s-} = D_tY_s$ is null with respect to
the marginal measure $\nu$.\\
 \item[c')] For $|\mu|$ almost all $(t_1,t_2)$,
\begin{equation}
    \lim_{\varepsilon \rightarrow 0} \frac{1}{\varepsilon} \int^{t_2 + \varepsilon}_{t_2} D_{t_1}Y_s ds \label{5.12'}
    \end{equation}
    exists a.s. It will be denoted by $(D_{t_1}Y_{t_2+})$. Moreover for each $s$, the set of $t$ such that $D_tY_{s+} = D_tY_s$ is null with respect to measure $\nu$.
\end{description}
If a), b), c) (resp. a), b), c')) are verified then $Y \in
Dom(\delta)$ and the forward integral $\int_0^T Y d^-X$ (resp. the
backward integral $\int_0^T Y d^+X$) exists and
\begin{equation}
\int_0^t Y d^-X = \int_0^t Y \delta X + \int_{[0,t]^2}  \label{52a}
D_{t_1}Y_{t_2-}d\mu(t_1,t_2)
\end{equation}
(resp.
\begin{equation}
\int_0^t Y d^+X = \int_0^t Y \delta X + \int_{[0,t]^2} \label{53a}
D_{t_1}Y_{t_2+}d\mu(t_1,t_2).)
\end{equation}
\end{thm}

\begin{rem} \label{r5.8}
i) Condition (\ref{5.11}) implies trivially
\begin{equation}
E\left( \int_{[0,T]^2} \left|D_{t_1}Y_{s_1}\right| \int_{[0,T]^2}
\left|D_{t_2}Y_{s_2}\right|\right) d\left| \mu
\right|(s_1,s_2)d\left| \mu \right|(t_1,t_2) < \infty. \label{5.13}
\end{equation}
ii) By Proposition \ref{P4} we know that $Y \in Dom(\delta)$.
\newline
iii) Taking into account the definition if $\int_0^t Y \delta X$, it will be enough to prove (\ref{52a}) and  (\ref{53a}) replacing $t$ with $T$.
\end{rem}

\begin{rem} \label{r8.3b}
In the case of Malliavin calculus for classical Brownian
motion, see Section \ref{sec3}, one has
\[
d\mu(s_1, s_2) = \delta(ds_2 - s_1) ds_1.
\]
So (\ref{5.13}) becomes
\begin{equation}
E\left( \iint_{[0,T]^2} \left| D_t Y_s \right|^2 ds dt \right) <
\infty. \label{5.14}
\end{equation}
\end{rem}

\vskip0.3cm

\begin{rem} \label{r8.4b}
Condition (\ref{5.12}) (resp. (\ref{5.12'}))  may be replaced by the existence a.s. of the trace $Tr Du$, where
\begin{equation}
\label{trace} Tr Du(t) =  \lim _{\varepsilon \to 0} \frac{1}{\varepsilon }\int
_{0}^{t} \langle DY_{s}, 1_{]s, s+\varepsilon ]} \rangle
_{{\cal{H}}}ds, \quad t \in [0,T].
\end{equation}
 This is a direct consequence of Fubini theorem. A similar condition related to symmetric integral appears in \cite{ALN}.
\end{rem}

\begin{lem} \label{l5.9}
Let $(Y_t)_{t \in [0,T]}$ be a process fulfilling points a), b), c) of
Theorem \ref{t5.7}. We set
\begin{equation}
Y^\varepsilon_t = \frac{1}{\varepsilon} \int_{(t - \varepsilon)^+}^t
Y_s ds \label{5.15}.
\end{equation}
Then  $Y^{\varepsilon}\in Dom(\delta)$ and for every $t$
\begin{equation}
\int_0^t Y^\varepsilon \delta X \stackrel{\varepsilon \rightarrow
0} {\longrightarrow} \int_0^t Y \delta X  \mbox{ in }
L^{2}(\Omega). \label{5.16}
\end{equation}

\end{lem}
{\bf Proof: } Again, in this proof it will be enough to set $t = T$. First, one can prove that if $Y \in Cyl(L_{\mu})$, $Y^{\varepsilon} \in Cyl(L_{\mu})$ and
\begin{equation}\label{E20}
D_tY_s^{\varepsilon} = \frac{1}{\varepsilon} \int_{s-\varepsilon} ^{s } D_t Y_r dr.
\end{equation}
Then we can establish that $Y^\varepsilon$ is a suitable limit of elements in $Cyl(L_\mu)$ so that $Y^{\varepsilon} \in |\mathbb{D}^{1,2}(L_\mu)|$.
We omit details at this level.
Relation (\ref{E20}) extends then to every $Y$ fulfilling the assumptions of the theorem.
According to Proposition \ref{P4}, $Y^\varepsilon \in Dom(\delta)$. Relation (\ref{E14}) in Proposition \ref{P4} gives
\begin{eqnarray} \label{E21}
E\left( \int_0^T (Y - Y^\varepsilon) \delta X\right)^2 \leq E(\|Y - Y^\varepsilon\|^2_\mathcal{H}) + E\left( \int_{[0,T]^2}
 d|\mu|(t_1,t_2)\| D_{\cdot}(Y_{t_1}-Y_{t_1}^\varepsilon) \|^2_{\mathcal{H}} \right).
\end{eqnarray}

We have to show that both expectations converge to zero.
The first expectation gives
\begin{equation}\label{E22}
E\left( \int_{[0,T]^2}d\mu(t_1,t_2) (Y_{t_1} - Y_{t_1}^\varepsilon)  (Y_{t_2} - Y_{t_2}^\varepsilon)\right).
\end{equation}
Using assumption a) of the theorem, Lebesgue dominated convergence theorem implies that (\ref{E22}) converges to
\[
E\left( \int_{[0,T]^2} d\mu(t_1,t_2) (Y_{t_1} - Y_{t_1-})  (Y_{t_2} - Y_{t_2-})\right).
\]
For each $\omega$ a.s. the discontinuities of $Y(\omega)$ are countable. The fact that $|\mu|$ is non-atomic implies that previous expectation is zero.
\par We discuss now the second expectation. It gives
\begin{equation} \label{E23}
\int_{[0,T]^2}d|\mu| (t_1,t_2) \int_{[0,T]^2} E \left( D_{s_1}(Y_{t_1} - Y_{t_1}^\varepsilon) D_{s_2} (Y_{t_1} - Y_{t_1}^\varepsilon)\right) d|\mu|(s_1,s_2).
\end{equation}
Taking into account assumptions b), c) of the theorem, previous term converges  to
\[
E \left( \int_{[0,T]^2}d|\mu| (t_1,t_2) \int_{[0,T]^2} d|\mu|(s_1,s_2) (D_{s_1}Y_{t_1}- D_{s_1}Y_{t_1-}) (D_{s_2} Y_{t_1} - D_{s_2}Y_{t_1-})\right).
\]
Using Cauchy-Schwarz this is bounded by
\[
E \left( \int_0^T d\nu (t) \int_0^T d\nu(s)(D_sY_t - D_sY_{t-})^2 \right).
\]
This quantity is zero because of c). \qed

\begin{rem} \label{r5.9}
If point c') is verified (instead of c) it is possible to state a similar version of the lemma with $Y^\varepsilon_t = \frac{1}{ \varepsilon }\int_t^{t+\varepsilon}Y_s ds$.
\end{rem}

It is interesting to observe that convergence (\ref{5.16}) holds weakly in $L^2(\Omega)$ even without assumption c). This constitutes the following proposition.
\begin{prop} \label{p5.9}
Let $(Y_t)_{t \in [0,T]}$ be a process fulfilling points a), b) of Theorem \ref{t5.7}. We set $Y^\varepsilon$ as in (\ref{5.15}). Then for every $t$,
\begin{equation}
\int_0^t Y^\varepsilon \delta X \stackrel{\varepsilon \rightarrow 0}{ \longrightarrow} \int_0^t Y \delta X \label{5.16bis}
\end{equation}
weakly in $L^2(\Omega)$.
\end{prop}
{\bf Proof: } We set again $t = T$. 
One can prove directly that $Y^\varepsilon$ belongs to $Dom(\delta)$ because of Fubini type Proposition \ref{p5.6}. Indeed, we set
\[
G = [0,T], \ \nu(ds) = ds, \ u(s,t) = Y_s 1_{]s,s+\varepsilon]}(t)
\]
and we verify the assumptions of the Proposition. Using Proposition \ref{P4} and points a), b) it is clear that $E(\int_0^T Y^\varepsilon \delta X)^2$ is bounded.
Then it is possible to show that the left term in
(\ref{5.16}) admits a subsequence $(\int_0^T
Y^{\varepsilon_n}\delta W)$ converging weakly to some square
integrable random variable $Z$.

\smallskip

Let $F \in Cyl$. By duality of Skorohod integral
\begin{eqnarray}
E\left(F \int_0^T Y^{\varepsilon_n} \delta X\right) &&= E \left(  \left\langle  DF,Y^{\varepsilon_n} \right\rangle_\mathcal{H} \right) \nonumber \\
&& = E \left(\int_{[0,T]^2} D_{s_1}F Y_{s_2}^{\varepsilon_n}
\mu(ds_1,ds_2)   \right) \nonumber \\
&& = E
\left(\int_{[0,T]^2} D_{s_1}F Y_{s_2-} \mu(ds_1,ds_2)   \right).
\nonumber
\end{eqnarray}
Now since $X $ is $L^{2}$ continuous, it is not difficult to see that the
\begin{equation}
|\mu|\left(\left\{s_1\right\} \times [0,T]\right) = |\mu| \left([0,T]
\times \{s_2\} \right) = 0.\label{5.17}
\end{equation}
Using Banach-Steinhaus theorem and the density of $Cyl$ in $L^2(\Omega)$, the convergence (\ref{5.16}) is established.
For $\omega$ a.s the set $N(\omega)$ of discontinuity of $Y(\omega)$
is countable. Consequently $|\mu|\left( [0,T] \times N(\omega)
\right) = 0$ and so
\begin{eqnarray}
E \left(\int_{[0,T]^2} D_{s_1}F Y_{s_2-} \mu(ds_1,ds_2)   \right)&& = E \left(\int_{[0,T]^2} D_{s_1}F Y_{s_2} \mu(ds_1,ds_2)   \right)\nonumber \\
&& =E\left(<DF,Y>_\mathcal{H} \right) = E\left(F \int_0^T Y \delta X
\right). \nonumber
\end{eqnarray}
\qed
\vskip0.5cm
{\bf Proof} of the Theorem \ref{t5.7}:  We set again $t = T$. We operate only for the forward integral. The backward case can be treated similarly.

Proposition \ref{p5.5} implies that
\begin{eqnarray}
I^-(\varepsilon, Y, dX, T) &&= \frac{1}{\varepsilon} \int_0^T ds Y_s \int_0^T 1_{]s, s+ \varepsilon]}(t) \delta X_t \nonumber\\
&& = \frac{1}{\varepsilon} \int_0^T ds \int_0^T Y_s  1_{]s, s+
\varepsilon]}(t) \delta X_t
+ \frac{1}{\varepsilon} \int_0^T ds \left( \int_{[0,T]^2} d\mu(t_1,t_2) D_{t_1}Y_s 1_{]s, s+ \varepsilon]}(t_2)\right)\nonumber\\
&&= I_1(T, \varepsilon) + I_2 (T, \varepsilon). \nonumber
\end{eqnarray}
  Proposition \ref{p5.6} says that
\[
I_1(T,\varepsilon) = \int_0^T Y_t^{\varepsilon} \delta X_t.
\]
According to Lemma \ref{l5.9}, $I_1(T, \varepsilon)$ converges in
$L^2(\Omega)$ to $\int_0^T Y \delta X$.

 We observe now that $I_2(T, \varepsilon)$ gives
\begin{equation}
\int_{[0,T]^2} d\mu(t_1,t_2) \frac{1}{\varepsilon} \int_{t_2 -
\varepsilon}^{t_2} ds D_{t_1} Y_s. \label{5.18}
\end{equation}
Assumptions b) and c) together with Lebesgue dominated convergence
theorem show that (\ref{5.18}) converges in $L^2(\Omega)$  to
\[
\int_{[0,T]^2} d\mu(t_1,t_2)D_{t_1}Y_{t_2}.
\]
\qed

\vskip0.5cm

In particular, we retrieve the result in Lemma \ref{wfb}.
\begin{cor} \label{c5.10}
Let $h$ be a cadlag function $h:[0,T] \rightarrow \mathbb{R}$. Then
\[
\int_0^T h dX = \int_0^T h d^-X = \int_0^T h d^+X.
\]
\end{cor}
{\bf Proof: } This is obvious because the Malliavin derivative of
$h$ vanishes. \qed.

\vskip0.5cm

\begin{cor} \label{c5.11}
Let $(Y_t)_{t  \in [0,T]}$ to be a process fulfilling assumptions
a), b), c), c') of Theorem {\ref{t5.7}}. Then the symmetric integral
of $Y$ with respect to $X$ is defined and
\begin{equation*}
\int_0^t Y d^oX  = \int_0^t Y \delta X + \frac{1}{2} \int_{[0,t]^2}
 \left( D_{t_1}  Y_{t_2+} + D_{t_1}Y_{t_2-} \right)d\mu(t_1,t_2)
\end{equation*}
 and
 \begin{equation*} [X,Y]_t = \int_{[0,t]^2}  \left( D_{t_1}
Y_{t_2+} + D_{t_1}Y_{t_2-} \right) d\mu(t_1,t_2).
\end{equation*}
\end{cor}

\begin{ex} \label{e5.12} \normalfont {\bf The case of a Gaussian martingale
$X$.}
\end{ex}
We recall by Section \ref{sec3} that $[X] = \lambda$, where $\lambda$ is a
deterministic increasing function vanishing at zero. Under
assumption a), b), c) of Theorem \ref{t5.7}
\[
\int_0^T Y d^-X = \int_0^T Y \delta X + \int_0^T d \lambda(t_1)
D_{t_1}Y_{t_1-}.
\]
Let $Y$ be $\mathcal{F}$-progressively measurable cadlag, such
that $\int_0^T Y_s^2 d[X]_s < \infty$ a.s. In \cite{RVsem} it is
also shown that $\int_0^T Y d^-X$ equals the It\^o integral
$\int_0^T Y dX$.
\par It is possible to see that the It\^o integral $Z = \int_0^T Y dX$ verifies the duality relation (\ref{5.5}) and so $Y \in Dom(\delta)$. Moreover
\[
\int_0^T Y d^{-}X = \int_0^T Y \delta X.
\]

\vskip0.5cm

 We will discuss now It\^o formula. Theorem \ref{t5.7} allows to state the following preliminary formulation.

\begin{lem} \label{p5.13}
Let $f \in C^2(\mathbb{R})$ such that $f''$ is bounded. Then, for every $t$,  $f'(X_t) \in Dom(\delta)$, and
\[
f(X_t) = f(X_0) + \int_0^t f'(X_s)\delta X_s + \int_{\Delta_t}
f''(X_{s_2})d\mu(s_1,s_2) + \frac{1}{2} \int_0^t f''(X_s)
d{\cal{E}}(s),
\]
where
\begin{equation*}
{\cal{E}}(t) =\mu(D_t), \hskip0.2cm  D_t =\left\{ (s,s)| s \leq t
\right\}, \hskip0.2cm  \Delta_t  =  \left\{ (s_1,s_2) | 0 \le  s_1 < s_2 \le t
\right\}.
\end{equation*}
\end{lem}
{\bf Proof: } It\^o formula for finite quadratic variation processes
was established for instance by \cite{RV2000}. It says
\[
f(X_t) = f(X_0) + \int_0^t f'(X) d^-X + \frac{1}{2} \int_0^t f''(X)
d[X].
\]
Now we need to apply Theorem \ref{t5.7}. For this we need to
verify its hypotheses. The assumption  a) is verified because
\[
\sup_{s \leq t} \left| f(X_s) \right| \leq \sup|f'|\sup_{t \leq
t}|X_s|.
\]

Since $X$ is a Gaussian process, (\ref{EsupG}) recalls that $\sup_{t
\leq T}|X_t| \in L^2(\Omega) $. On the other hand, setting $Y_t = f'(X_t)$,
\[D_{t_1}Y_{t_2} = f''(X_{t_2}) 1_{]0,t_2[}(t_1)
\]
and so b) is also verified.
\[
\lim_{\varepsilon \rightarrow 0} \int_{t_2 - \varepsilon}^{t_2}
D_{t_1}Y_s ds = \left\{
\begin{array}{clrr}
0 &\ \ \  t_1 \geq t_2 \\
f''(X_{t_2})& \ \ \ t_1 < t_2
\end{array}
\right.
\]
and c) is verified. Therefore
\[
\int_0^t f'(X) d^- X = \int_0^t f'(X) \delta X + \int_{\Delta_t}
f''(X_{t_2})d\mu(t_1,t_2).
\]
Moreover, by Lemmas \ref{l2.6} and \ref{l2.7} we have
\[
\frac{1}{2} \int_0^t f''(X_s)d[X]_s = \frac{1}{2} \int_0^t
f''(X_{s}) d{\cal{E}}(s).
\]\qed

\par The statement of Lemma \ref{p5.13} can be precised better.

\begin{lem}\label{lemmaA}
Let $\gamma (t) = R(t,t)$. For every $\psi \in C^0(\mathbb{R}_+)$,
\begin {equation}
\int_0^t \psi d^- \gamma = \int_0^t \psi(s) d\mathcal{E}(s)+ 2 \int_{\Delta_t} \psi(s_2)d\mu(s,s_2), \label{eq1}
\end{equation}
where
\begin{equation}
\int_0^t \psi d^- \gamma = \lim_{\varepsilon \rightarrow 0} \int _0^t \psi (s) \frac{\gamma(s+\varepsilon)+ \gamma(s)}{\varepsilon} ds \label{eq2}
\end{equation}
pointwise. Moreover $\gamma$ is a bounded variation function. 
\end{lem}
{\bf Proof:  } The integral inside the right-hand side of (\ref{eq2}) gives
\[
\int_0^t \psi(s) \frac{R(s+\varepsilon,s+\varepsilon)-R(s,s)}{\varepsilon}ds = (I_1 + 2 I_2)(\varepsilon,t),
\]
where
\begin{eqnarray}
I_1(\varepsilon,t)&=& \int_0^t \psi(s) \Delta_{]s,s+\varepsilon]^2} R \frac{ds}{\varepsilon}, \nonumber\\
I_2(\varepsilon,t) &=&  \int_0^t \frac{R(s+\varepsilon,s)-R(s,s)}{\varepsilon}\psi(s)ds. \nonumber
\end{eqnarray}
We have
\[
I_1(\varepsilon,t) = \int_0^t \psi(s) E(X_{s+\varepsilon}-X_s)^2 \frac{ds}{\varepsilon}= \int_0^t \psi(s) d(E(C_\varepsilon(X,X,s)))ds.
\]
Taking into account Definition \ref{d2.6} and the fact that $\mathcal{E}$ is a continuous function,
\[
I_1(\varepsilon,t) \longrightarrow_{\varepsilon \rightarrow 0} \int_0^t \psi(s) d\mathcal{E}(s), \]
for any $t \in [0,T]$.
It remains to control $I_2(\varepsilon,t)$.
\par Since 
\[
R(s+\varepsilon,s)-R(s,s) = \int_{[0,t]^2} 1_{]0,s] \times ]s,s+\varepsilon]}d\mu(s_1,s_2),
\]
we have
\[
I_2(\varepsilon,t)= \int_{[0,t]^2} d\mu(s_1,s_2) \int_{]s_1 \vee (s_2-\varepsilon)^+, s_2]} \psi(s) \frac{ds}{\varepsilon},
\]
with the convention that $]a,b] = \varnothing$ if $b\leq a$.
\par We distinguish two cases:
\begin{itemize}
	\item If $s_2 \leq s_1$, then
	\[
	I_2(\varepsilon,t) = 0.
	\]
	\item If $s_2 > s_1$, then
	\[
	\frac{1}{\varepsilon}\int_{]s_1 \vee (s_2-\varepsilon),s_2]}\psi(s)ds \longrightarrow_{\varepsilon \rightarrow 0} \psi(s_2).
	\]pointwise.
	\end{itemize}
\par
Finally
\[
I_2(\varepsilon,t) \longrightarrow \int_{\Delta_t}d\mu(s_1,s_2)\psi(s_2).
\]
and the proof of (\ref{eq1}) is established.
\par In particular
\[
\gamma_t = \mathcal{E}(t) + 2 \mu(\Delta_t).
\]
Since $\mathcal{E}$ has bounded variation and the total variation of $t \longmapsto \mu(\Delta_t)$ is bounded by $|\mu|([0,T]^2)$, then $\gamma$ is also a bounded variation function.     \qed

\begin{cor}\label{cor8.11c}
Let $f \in C^2(\mathbb{R})$ such that $f''$ is bounded. We set $\gamma_t = Var(X_t)$. Then 
\[
f(X_t)=f(X_0)+\int_0^t f'(X_s)\d X_s + \frac{1}{2}\int_0^t f''(X_s)d\gamma(s).
\]
\end{cor}
{\bf Proof:   } It follows from Lemma \ref{p5.13} and Lemma \ref{lemmaA} setting $\psi(t)=f''(X_t)$. \qed

\vskip0.5cm

 We would like to examine some particular cases. For this
we decompose $\mu$ into $\mu_{d} + \mu_{od}$ where for every $A \in
B([0,T]^2)$

\begin{equation*}
\mu_d(A) = \mu(A \cap D_T), \hskip 0.3cm  \mu_{od}(A) = \mu(A
\backslash D_T).
\end{equation*}
Hence $\mu_d$ is concentrated on the diagonal, $\mu_{od}$ outside
the diagonal.

\smallskip

 We recall that
\[
{\cal{E}}(t) = \mu(D_t),
\]
where ${\cal{E}} = {\cal{E}}(X)$ is the energy function defined in
Section \ref{sec2}. Consider the repartition functions $R_d$,
$R_{od}$ of $\mu_d$, $\mu_{od}$. We have
\begin{equation*}
 R_d(s_1,s_2) = \mu_d\left(   ]0,s_1] \times ]0,s_2] \right) = {\cal{E}}(s_1 \wedge s_2)
 \end{equation*}
 and
 \begin{equation*}
  R_{od}(s_1,s_2) = \mu_{od}\left(   ]0,s_1] \times ]0,s_2]
\right).  \end{equation*}

\vskip0.5cm

\begin{rem} \label{r5.14}
i) Setting $\psi = {\cal{E}}$, there is a Gaussian martingale $M$  with
covariance $R_d$.
It is enough to take $M_t = W_{\psi(t)},$ where $(W_t)$ is a classical Brownian motion.\\
ii) $R_{od}(s_1,s_2) = Cov(X_{s_1},X_{s_2}) - Cov(M_{s_1},M_{s_2}).$
\end{rem}

\vskip0.5cm

\begin{prop} \label{p5.15}
Suppose that
\begin{description}
    \item[i)] ${\cal{E}}$ is absolutely continuous, i.e. there is a locally integrable function  ${\cal{E}}' $  such that
    \[
    {\cal{E}}(t) = \int_0^t {\cal{E}}'(s)ds:
    \]
    \item[ii)] $\mu_{od}$ is absolutely continuous with respect to Lebesgue measure. In particular one has
    \[
    \mu_{od}(]0,s_1] \times]0,s_2]) = \int_{]0,s_1] \times]0,s_2]}\frac{\partial^2 R}{\partial s_1 \partial s_2} ds_1 ds_2.
    \]
    In fact it is clear that
    \[
    \frac{\partial^2 R}{\partial s_1 \partial s_2} = \frac{\partial^2 R_{od}}{\partial s_1 \partial     s_2} \textrm{ on } [0,T]^2 \backslash D_T.
    \]
\end{description}
Then the conclusion of the Theorem \ref{t5.7} holds replacing
assumptions
\begin{itemize}
    \item Y is cadlag,
    \item c) (resp. c'),
\end{itemize}
with
\begin{itemize}
    \item For t a.e. Lebesgue
    \begin{equation} \label{5.19}
        D_t Y_{t-} = \lim_{\varepsilon \rightarrow 0} \frac{1}{\varepsilon} \int_{t - \varepsilon}^{t} D_t Y_s ds \textrm{ exists a.s.}
    \end{equation}
    (resp.
        \begin{equation} \label{5.19'}
        D_t Y_{t+} = \lim_{\varepsilon \rightarrow 0} \frac{1}{\varepsilon} \int_{t}^{t + \varepsilon} D_t Y_s ds \textrm{ exists a.s.})
    \end{equation}
\end{itemize}
Moreover the conclusion of the theorem can be stated as
\begin{equation} \label{5.20}
\int_0^T Y d^-X = \int_0^T Y \delta X + \int_0^T D_tY_{t-}
{\cal{E}}'(t)dt + \int_{[0,T]^2} D_{t_1}Y_{t_2} \frac{\partial^2
R}{\partial s_1
\partial     s_2}(t_1,t_2)dt_1,dt_2.
\end{equation}
(resp.
\begin{equation} \label{5.20'}
\int_0^T Y d^+X = \int_0^T Y \delta X +\int_0^T D_tY_{t+}
{\cal{E}}'(t)dt  +\int_{[0,T]^2} D_{t_1}Y_{t_2} \frac{\partial^2
R}{\partial s_1
\partial     s_2}(t_1,t_2)dt_1,dt_2.)
\end{equation}
\end{prop}

\vskip0.5cm

\begin{rem} \label{r5.16}
If c) and c') of Theorem \ref{t5.7}, with (\ref{5.19}) and (\ref{5.19'}) are
verified then
\begin{equation*}
\int_0^t Y d^oX = \int_0^t Y \delta X + \int_0^t (D_sY_{s+} +
D_sY_{s-}) {\cal{E}}'(s)dt + \int_{[0,t]^2} D_{t_1}Y_{t_2}
\frac{\partial^2 R}{\partial s_1 \partial s_2} (t_1,t_2)dt_1dt_2
\end{equation*}
and \begin{equation*} [X,Y]_{t} = \int_0^t (D_sY_{s+} + D_sY_{s-})
{\cal{E}}'(s)dt.
\end{equation*}
\end{rem}

\vskip0.5cm

 Before proceeding to the proof, we recall a basic result
which can be found in \cite{Stein}.
\begin{lem} \label{l5.17}
Let $g \in L^p(\mathbb{R})$, $1 \leq p < \infty$. We set
\[
g_{\varepsilon}(x) = \frac{1}{\varepsilon} \int_{x - \varepsilon}^x
g(y) dy \mbox { or } g_{\varepsilon}(x) = \frac{1}{\varepsilon}
\int^{x + \varepsilon}_x g(y) dy.
\]
Then $g_{\varepsilon} \rightarrow g$ a.e. and in $L^p$.
\end{lem}
{\bf Proof} of Proposition \ref{p5.15}: The proof follows the same
line as the proof of Theorem \ref{t5.7}. \\
a) First we need to adapt
Lemma \ref{l5.9} to show that $\lim_{\varepsilon\rightarrow 0}
\int_0^T Y_s^\varepsilon \delta X_s = \int_0^T Y_s \delta X_s$ in
$L^2(\Omega)$ where $Y^{\varepsilon}$ still denotes the same
approximation process. Again we need to show that the right hand
side of  (\ref{E21}) converges to zero when $\varepsilon \rightarrow
0$. Its first term  gives $(I_1 + I_2)(\varepsilon)$, where
\begin{eqnarray}
I_1(\varepsilon) &=& E \left(\int_{[0,T]^2} ds_1 ds_2 \frac{\partial^2R }{\partial s_1
 \partial s_2 } (Y_{s_1}^{\varepsilon} -Y_{s_1}) (Y_{s_2}^{\varepsilon} -Y_{s_2}) \right), \nonumber\\
I_2(\varepsilon) &=& E \left(\int_0^T ds \mathcal{E}'(s) (Y_{s}^{\varepsilon} -Y_s)^2 ) \right) \nonumber.
\end{eqnarray}
Lemma \ref{l5.17} implies that $Y^\varepsilon \longrightarrow Y$
a.e. $dP \otimes Leb$. Lebesgue dominated convergence theorem and
Assumption a) imply that $I_1(\varepsilon) \longrightarrow 0$ and
$I_2(\varepsilon) \longrightarrow 0$, when $\varepsilon$
converges to zero. It remains to control the second term in
(\ref{E21}) which is given by (\ref{E23}). This second term  gives
$K_1(\varepsilon) + K_2(\varepsilon)$ with
\begin{eqnarray}
K_1(\varepsilon) &=& E \left(\int_{[0,T]^2} ds_1 ds_2 \left| \frac{\partial^2R }{\partial s_1 \partial s_2 } \right| \int_{[0,T]^2} dt_1 dt_2 \left|
\frac{\partial^2R }{\partial t_1 \partial t_2 } \right| |D_{s_1}Y_{t_1}^{\varepsilon} -D_{s_1}Y_{t_1}| |D_{s_2}Y_{t_1}^{\varepsilon} -D_{s_2}Y_{t_1}| \right), \nonumber\\
K_2(\varepsilon) &=& E \left(\int_0^T ds \mathcal |{\cal E}'(s)| \int_0^T dt (D_sY_{t}^{\varepsilon} -D_sY_t)^2 ) \right) \nonumber.
\end{eqnarray}
Point b) and (\ref{5.19}) allow to show that $K_1(\varepsilon) + K_2(\varepsilon) \longrightarrow 0$.

b) The other point concerns the convergence of $I_2(T,
\varepsilon)$ appearing in the proof of Theorem 8.1. To prove the
convergence of (\ref{5.18}) we separate again $\mu = \mu_d +
\mu_{od}$ and we use (\ref{5.19}) on the diagonal. Finally Lemma
\ref{l5.17}, (\ref{5.11}) and Lebesgue dominated convergence theorem
show that for $t_1, \; t_2$, $t_1 \neq t_2$ a.e.
\[
\frac{1}{\varepsilon} \int_{t_2}^{t_2 + \varepsilon} ds D_{t_1}Y_s
\stackrel{\varepsilon \rightarrow 0}{\longrightarrow}
 D_{t_1} Y_{t_2} \ \textrm{ a.e. }\  dP \otimes dt_1 dt_2.\] \qed

 \vskip0.5cm

\begin{ex} Let us apply the  obtained results to some particular
examples.
\end{ex}
\begin{description}
    \item[a)] {\em Case of a Gaussian martingale with absolutely continuous quadratic variation
    $\lambda(t) = \lambda(0) + \int_0^t \dot{\lambda}(s) ds$.}
    \begin{eqnarray}
    &&R(t_1,t_2)  =  \lambda(t_1 \wedge t_2),\ \ \gamma(t) = Var(X_t) = \lambda(t), \nonumber\\
    &&{\cal{E}}'(t) = \dot{\lambda}(t), \ \ R = R_d,\nonumber\\
    &&\frac{\partial ^2 R}{\partial t_1 \partial t_2}  =  0 \textrm{  a.e. Lebesgue.} \nonumber
    \end{eqnarray}
    Let $Y$ be as in Proposition \ref{p5.15}. Then
    \[
    \int_0^T Y d^- X = \int_0^T Y \delta X + \int_0^T  D_tY_{t-} \dot{\lambda}(t) dt.
    \]

    \item[b)] {\em The case of fractional Brownian motion $H > 1/2$.}
    \\
    We have
    \begin{eqnarray}
    && R = R_{od}\nonumber\\
    && \frac{\partial ^2 R}{\partial t_1 \partial t_2} = 2H (2H - 1)|t_2 -t_1|^{2H-2}, \ \gamma(t)=t^{2H}.
    \end{eqnarray}
    One obtains the classical results for fractional Brownian motion as in \cite{AN}, for instance
    \[
    \int_0^T Y d^-X = \int_0^T Y \delta X + H(2H - 1) \int_{[0,T]^2} D_{t_1}Y_{t_2}\vert t_2 - t_1\vert
    ^{2H -2}dt_1dt_2.
    \]
    Corollary \ref{cor8.11c} provides the following It\^o formula:
    \begin{equation}\label{ito-fbm}
    f(X_t) = f(X_0) + \int_0^t f'(X) \delta X + H\int_0^t  f''(X_s)s^{2H - 1}ds.
    \end{equation}

    \item[c)] {\em The case of bifractional Brownian motion $X = B^{H,K}$, $HK \geq \frac{1}{2}$. } \\
    It is easy to verify that $Var(X_t) = t^{2HK}$ so that
    \par Corollary \ref{cor8.11c} implies 
\[
f(X_{t})= f(0)+\int_{0}^{t}f'(X_{s})\delta X_{s} + HK\int _{0}^{t}
f''(X_{s})s^{2HK-1}ds.
\]
In particular, if $HK=\frac{1}{2}$ we get a formula which looks very similar to the one for the Brownian motion:
\[
f(X_t)  = f(X_0) + \int_0^t f'(X_s) \delta X_s +\frac{1}{2} \int_0^t
    f''(X_s)ds.
\]
\end{description}

\end {document}